# AN INVARIANCE PRINCIPLE FOR SEMIMARTINGALE REFLECTING BROWNIAN MOTIONS IN DOMAINS WITH PIECEWISE SMOOTH BOUNDARIES[1]

BY W. KANG AND R. J. WILLIAMS

*Carnegie Mellon University and University of California, San Diego*

Semimartingale reflecting Brownian motions (SRBMs) living in the closures of domains with piecewise smooth boundaries are of interest in applied probability because of their role as heavy traffic approximations for some stochastic networks. In this paper, assuming certain conditions on the domains and directions of reflection, a perturbation result, or *invariance principle,* for SRBMs is proved. This provides sufficient conditions for a process that satisfies the definition of an SRBM, except for small random perturbations in the defining conditions, to be close in distribution to an SRBM. A crucial ingredient in the proof of this result is an oscillation inequality for solutions of a perturbed Skorokhod problem. We use the invariance principle to show weak existence of SRBMs under mild conditions. We also use the invariance principle, in conjunction with known uniqueness results for SRBMs, to give some sufficient conditions for validating approximations involving (i) SRBMs in convex polyhedrons with a constant reflection vector field on each face of the polyhedron, and (ii) SRBMs in bounded domains with piecewise smooth boundaries and possibly nonconstant reflection vector fields on the boundary surfaces.

**1. Introduction.** Semimartingale reflecting Brownian motions (SRBMs) living in the closures of domains with piecewise smooth boundaries are of interest in applied probability because of their role as heavy traffic diffusion approximations for some stochastic networks. The nonsmoothness of the boundary for such a domain, combined with discontinuities in the oblique directions of reflection at intersections of smooth boundary surfaces, present

Received May 2006; revised November 2006.
[1]Supported in part by NSF Grants DMS-03-05272 and DMS-06-04537.
*AMS 2000 subject classifications.* 60F17, 60J60, 60K25, 90B15, 93E03.
*Key words and phrases.* Semimartingale reflecting Brownian motion, piecewise smooth domain, invariance principle, oscillation inequality, Skorokhod problem, stochastic networks.

challenges in the development of a rigorous theory of existence, uniqueness and approximation for such SRBMs.

When the state space is an orthant and the direction of reflection is constant on each boundary face, a necessary and sufficient condition for weak existence and uniqueness of an SRBM is known [14]. This condition involves a so-called completely-$\mathcal{S}$ condition on the matrix formed by the reflection directions. An invariance principle for such SRBMs was established in [15] and used in [16] to justify heavy traffic diffusion approximations for certain open multiclass queueing networks. Loosely speaking, the invariance principle shows that, assuming uniqueness in law for the SRBM, a process satisfying the definition of the SRBM, except for small perturbations in the defining conditions, is close in distribution to the SRBM.

For more general domains with piecewise smooth boundaries, some conditions for existence and uniqueness of SRBMs are known. In particular, for convex polyhedrons with a constant direction of reflection on each boundary face, necessary and sufficient conditions for weak existence and uniqueness of SRBMs are known for simple convex polyhedrons (where precisely $d$ faces meet at each vertex in $d$-dimensions) and sufficient conditions are known for nonsimple convex polyhedrons, see [4]. For a bounded domain that can be represented as a finite intersection of domains, each of which has a $C^1$-boundary and an associated uniformly Lipschitz continuous reflection vector field, sufficient conditions for strong existence and uniqueness were provided by Dupuis and Ishii [6]; in fact, these authors study stochastic differential equations with reflection which include SRBMs. Despite these existence and uniqueness results, a general invariance principle for SRBMs living in the closures of domains with piecewise smooth boundaries has not been proved to date. (We note that for the special case when the directions of reflection are normal, that is, perpendicular to the boundary, there are a number of perturbation results for reflecting Brownian motions. Our emphasis here is on treating a wide range of oblique reflection directions.)

Motivated by its potential for use in approximating heavily loaded stochastic networks that are more general than open multiclass queueing networks, in this paper, we formulate and prove an invariance principle for SRBMs living in the closures of domains with piecewise smooth boundaries with possibly nonconstant directions of reflection on each of the smooth boundary surfaces. An application of the results of this paper to the analysis of an internet congestion control model can be found in [13]. An outline of the current paper is as follows.

The definition of an SRBM and assumptions on the domains and directions of reflection are given in Sections 2 and 3, respectively. Some sufficient conditions for these assumptions to hold are provided in Section 3. Section 4 is devoted to proving the main result of this paper, namely, the *invariance principle*. A key element for our proof of this result is an oscillation inequality for solutions of a perturbed Skorokhod problem; this inequality, which may be of independent interest, is proved in Section 4.1. In Section 5 we give some applications of the invariance principle. We prove weak existence of SRBMs under the conditions specified in Section 3. We also use the invariance principle, in conjunction with known uniqueness results for SRBMs, to give sufficient conditions for validating approximations involving (i) SRBMs in convex polyhedrons with a constant reflection vector field on each face of the polyhedron, and (ii) SRBMs in bounded domains with piecewise smooth boundaries and possibly nonconstant reflection vector fields on the boundary surfaces.

Beyond its possible use in justifying SRBM approximations for stochastic networks, the invariance principle might be used to justify numerical approximations to SRBMs. A further possible extension of the results stated here would involve an invariance principle for stochastic differential equations with reflection. The oscillation inequality for the perturbed Skorokhod problem and associated criteria for $C$-tightness described in Sections 4.1 and 4.2 are likely to be useful for this. We have not developed such an extension here as that would involve introduction of extra assumptions that would make the result less relevant for potential applications to stochastic networks. In particular, the approximating processes would involve stochastic integrals driven by a Brownian motion, whereas in stochastic network applications, the Brownian motion typically only appears in the limit.

1.1. *Notation, terminology and preliminaries.* Let $\mathbb{N}$ denote the set of all positive integers, that is, $\mathbb{N} = \{1, 2, \ldots\}$, $\mathbb{R}$ denote the set of real numbers, which is also denoted by $(-\infty, \infty)$, $\mathbb{R}_+$ denote the nonnegative half-line, which is also denoted by $[0, \infty)$. For $x \in \mathbb{R}$, we write $|x|$ for the absolute value of $x$, $[x]$ for the largest integer less than or equal to $x$, $x^+$ for the positive part of $x$. For any positive integer $d$, we let $\mathbb{R}^d$ denote $d$-dimensional Euclidean space, where any element in $\mathbb{R}^d$ is denoted by a column vector. Let $\|\cdot\|$ denote the Euclidean norm on $\mathbb{R}^d$, that is, $\|x\| = (\sum_{i=1}^d x_i^2)^{1/2}$ for $x \in \mathbb{R}^d$, and $\langle \cdot, \cdot \rangle$ denote the inner product on $\mathbb{R}^d$, that is, $\langle x, y \rangle = \sum_{i=1}^d x_i y_i$, for $x, y \in \mathbb{R}^d$. We note that for any $x \in \mathbb{R}^d$, $\|x\| \leq \sum_{i=1}^d |x_i|$. Let $\mathbb{R}_+^d$ denote the positive orthant in $\mathbb{R}^d$, that is, $\mathbb{R}_+^d = \{x \in \mathbb{R}^d : x_i \geq 0, 1 \leq i \leq d\}$. Let $\mathcal{B}(S)$ denote the Borel $\sigma$-algebra on $S \subset \mathbb{R}^d$, that is, the collection formed by intersecting all Borel sets in $\mathbb{R}^d$ with $S$. Let $\text{dist}(x, S)$ denote the distance between $x \in \mathbb{R}^d$ and $S \subset \mathbb{R}^d$, that is, $\text{dist}(x, S) = \inf\{\|x - y\| : y \in S\}$, with the convention that $\text{dist}(x, \varnothing) = \infty$ for $x \in \mathbb{R}^d$. Let $U_r(S)$ denote the closed set $\{x \in \mathbb{R}^d : \text{dist}(x, S) \leq r\}$ for any $r > 0$ and $S \subset \mathbb{R}^d$, where if $S = \varnothing$, $U_r(S) = \varnothing$ for all $r > 0$. Let $B_r(x)$ denote the closed ball $\{y \in \mathbb{R}^d : \|y - x\| \leq r\}$ for any $x \in \mathbb{R}^d$ and $r > 0$. For any set $S \subset \mathbb{R}^d$, we write $\overline{S}$ for the closure of $S$, $S^o$ for the interior of $S$ and $\partial S = \overline{S} \setminus S^o$. For a finite set $S$, $|S|$ denotes the number of elements in $S$. For any $v \in \mathbb{R}^d$, $v'$ denotes the transpose of $v$. Inequalities between vectors in $\mathbb{R}^d$ should be interpreted componentwise, that is, if $u, v \in \mathbb{R}^d$, then $u \leq (<) v$ means that $u_i \leq (<) v_i$ for each $i \in \{1, \ldots, d\}$. For any matrix $A$, let $A'$ denote the transpose of $A$. For any function $x : \mathbb{R}_+ \to \mathbb{R}^d$, $x(t-)$ denotes the left limit of $x$ at $t > 0$ whenever $x$ has a left limit at $t$; unless explicitly stated otherwise, $x(0-) \equiv 0$, where $0$ is the zero vector in $\mathbb{R}^d$. For any function $x : \mathbb{R}_+ \to \mathbb{R}^d$, we let $\Delta x(t) = x(t) - x(t-)$ for $t \in \mathbb{R}_+$ when $x(t-)$ exists. We let $\mathbf{0}$ be the constant deterministic function $x : \mathbb{R}_+ \to \mathbb{R}^d$ such that $x(t) = 0$ for all $t \in \mathbb{R}_+$.

A domain in $\mathbb{R}^d$ is an open connected subset of $\mathbb{R}^d$. For each continuously differentiable function $f$ defined on some nonempty domain $S \subset \mathbb{R}^d$, $\nabla f(x)$ is the gradient of $f$ at $x \in S$. For each $x \in \mathbb{R}^d$, a neighborhood $V_x$ of $x$ is a bounded domain in $\mathbb{R}^d$ that contains $x$. For any nonempty domain $S \subset \mathbb{R}^d$, we say that *the boundary $\partial S$ of $S$ is $C^1$*, if for each $x \in \partial S$ there exists a Euclidean coordinate system $C_x$ for $\mathbb{R}^d$ centered at $x$, an $r_x > 0$, and a once continuously differentiable function $\varphi_x : \mathbb{R}^{d-1} \to \mathbb{R}$ such that $\varphi_x(0) = 0$ and

$$S \cap B_{r_x}(x) = \{z = (z_1, \ldots, z_d)' \text{ in } C_x : z_d > \varphi_x(z_1, \ldots, z_{d-1})\} \cap B_{r_x}(x).$$

Then, for $x \in \partial S$, the inward unit normal to $\partial S$ at $z \in \partial S \cap B_{r_x}(x)$ is given in the coordinate system $C_x$ by

$$n(z) = \frac{1}{(1 + \|\nabla \varphi_x(z_1, \ldots, z_{d-1})\|^2)^{1/2}} (-\nabla \varphi_x(z_1, \ldots, z_{d-1})', 1)',$$

where $\nabla \varphi_x(z_1, \ldots, z_{d-1}) = (\frac{\partial \varphi_x}{\partial z_1}, \ldots, \frac{\partial \varphi_x}{\partial z_{d-1}})'(z_1, \ldots, z_{d-1})$. For any nonempty convex set $S \subset \mathbb{R}^d$, we call a vector $n \in \mathbb{R}^d \setminus \{0\}$ an inward unit normal vector to $S$ at $x \in \partial S$ if $\|n\| = 1$ and $\langle n, y - x \rangle \geq 0$ for all $y \in S$. Note that such a vector need not be unique.

All stochastic processes used in this paper will be assumed to have paths that are right continuous with finite left limits (abbreviated henceforth as r.c.l.l.). A process is called continuous if almost surely its sample paths are continuous. We denote by $D([0, \infty), \mathbb{R}^d)$ the space of r.c.l.l. functions from $[0, \infty)$ into $\mathbb{R}^d$ and we endow this space with the usual Skorokhod $J_1$-topology (cf. Chapter 3 of [7]). We denote by $C([0, \infty), \mathbb{R}^d)$ the space of continuous functions from $[0, \infty)$ into $\mathbb{R}^d$. The Borel $\sigma$-algebra on either $D([0, \infty), \mathbb{R}^d)$ or $C([0, \infty), \mathbb{R}^d)$ will be denoted by $\mathcal{M}^d$. The abbreviation *u.o.c.* will stand for *uniformly on compacts* and will be used to indicate that a sequence of functions in $D([0, \infty), \mathbb{R}^d)$ (or $C([0, \infty), \mathbb{R}^d)$) is converging uniformly on compact time intervals to a limit in $D([0, \infty), \mathbb{R}^d)$ (or $C([0, \infty), \mathbb{R}^d)$). Consider $W^1, W^2, \ldots, W$, each of which is a $d$-dimensional process (possibly defined on different probability spaces). The sequence $\{W^n\}_{n=1}^\infty$ is said to be *tight* if the probability measures induced by the $W^n$ on the measurable space $(D([0, \infty), \mathbb{R}^d), \mathcal{M}^d)$ form a tight sequence, that is, they form a weakly relatively compact sequence in the space of probability measures on $(D([0,\infty), \mathbb{R}^d), \mathcal{M}^d)$. The notation "$W^n \Rightarrow W$" will mean that, as $n \to \infty$, the sequence of probability measures induced on $(D([0,\infty), \mathbb{R}^d), \mathcal{M}^d)$ by $\{W^n\}$ converges weakly to the probability measure induced on the same space by $W$. We shall describe this in words by saying that $W^n$ converges weakly (or in distribution) to $W$ as $n \to \infty$. The sequence of processes $\{W^n\}_{n=1}^\infty$ is called *C-tight* if it is tight, and if each weak limit point, obtained as a weak limit along a subsequence, almost surely has sample paths in $C([0,\infty), \mathbb{R}^d)$. The following proposition provides a useful criterion for checking $C$-tightness.

PROPOSITION 1.1. *Suppose that, for each $n \in \mathbb{N}$, $W^n$ is a $d$-dimensional process defined on the probability space $(\Omega^n, \mathcal{F}^n, P^n)$. The sequence $\{W^n\}_{n=1}^\infty$ is C-tight if and only if the following two conditions hold:*

(i) *For each $\eta > 0$ and $T \geq 0$, there exists a finite constant $M_{\eta,T} > 0$ such that*

$$\liminf_{n \to \infty} P^n \left\{ \sup_{0 \leq t \leq T} \|W^n(t)\| \leq M_{\eta,T} \right\} \geq 1 - \eta. \tag{1}$$

(ii) *For each $\varepsilon > 0$, $\eta > 0$ and $T > 0$, there exists $\lambda \in (0, T)$ such that*

$$\limsup_{n \to \infty} P^n \{ w_T(W^n, \lambda) \geq \varepsilon \} \leq \eta, \tag{2}$$

*where for $x \in D([0,\infty), \mathbb{R}^d)$,*

$$w_T(x, \lambda) = \sup \left\{ \sup_{u,v \in [t, t+\lambda]} \|x(u) - x(v)\| : 0 \leq t < t + \lambda \leq T \right\}. \tag{3}$$

PROOF. See Proposition VI.3.26 in [12]. □

A $d$-dimensional process $W$ is said to be *locally of bounded variation* if all sample paths of $W$ are of bounded variation on each finite time interval. For such a process $W$, we define $\mathcal{V}(W) = \{\mathcal{V}(W)(t), t \geq 0\}$ such that for each $t \geq 0$,

$$\mathcal{V}(W)(t) = \|W(0)\|$$
$$+ \sup \left\{ \sum_{i=1}^{l} \|W(t_i) - W(t_{i-1})\| : 0 = t_0 < t_1 < \cdots < t_l = t, l \geq 1 \right\}.$$

A triple $(\Omega, \mathcal{F}, \{\mathcal{F}_t, t \geq 0\})$ will be called a filtered space if $\Omega$ is a set, $\mathcal{F}$ is a $\sigma$-algebra of subsets of $\Omega$, and $\{\mathcal{F}_t, t \geq 0\}$ is an increasing family of sub-$\sigma$-algebras of $\mathcal{F}$, that is, a filtration. Henceforth, the filtration $\{\mathcal{F}_t, t \geq 0\}$ will be simply written as $\{\mathcal{F}_t\}$. If $P$ is a probability measure on $(\Omega, \mathcal{F})$, then $(\Omega, \mathcal{F}, \{\mathcal{F}_t\}, P)$ is called a filtered probability space. A $d$-dimensional process $X = \{X(t), t \geq 0\}$ defined on $(\Omega, \mathcal{F}, P)$ is called $\{\mathcal{F}_t\}$-adapted if for each $t \geq 0$, $X(t): \Omega \to \mathbb{R}^d$ is measurable when $\Omega$ is endowed with the $\sigma$-algebra $\mathcal{F}_t$.

Given a filtered probability space $(\Omega, \mathcal{F}, \{\mathcal{F}_t\}, P)$, a vector $\mu \in \mathbb{R}^d$, a $d \times d$ symmetric, strictly positive definite matrix $\Gamma$, and a probability distribution $\nu$ on $(\mathbb{R}^d, \mathcal{B}(\mathbb{R}^d))$, an $\{\mathcal{F}_t\}$-Brownian motion with drift vector $\mu$, covariance matrix $\Gamma$, and initial distribution $\nu$, is a $d$-dimensional $\{\mathcal{F}_t\}$-adapted process defined on $(\Omega, \mathcal{F}, \{\mathcal{F}_t\}, P)$ such that the following hold under $P$:

(a) $X$ is a $d$-dimensional Brownian motion whose sample paths are almost surely continuous and that has initial distribution $\nu$,

(b) $\{X_i(t) - X_i(0) - \mu_i t, \mathcal{F}_t, t \geq 0\}$ is a martingale for $i = 1, \ldots, d$, and

(c) $\{(X_i(t) - X_i(0) - \mu_i t)(X_j(t) - X_j(0) - \mu_j t) - \Gamma_{ij} t, \mathcal{F}_t, t \geq 0\}$ is a martingale for $i, j = 1, \ldots, d$.

In this definition, the filtration $\{\mathcal{F}_t\}$ may be larger than the one generated by $X$; however, for each $t \geq 0$, under $P$, the $\sigma$-algebra $\mathcal{F}_t$ is independent of the increments of $X$ from $t$ onward. The latter follows from the martingale properties of $X$. If $\nu = \delta_x$, the unit mass at $x \in \mathbb{R}^d$, we say that $X$ starts from $x$.

**2. Definition of an SRBM.** Let $G = \bigcap_{i \in \mathcal{I}} G_i$ be a nonempty domain in $\mathbb{R}^d$, where $\mathcal{I}$ is a nonempty finite index set and for each $i \in \mathcal{I}$, $G_i$ is a nonempty domain in $\mathbb{R}^d$. For simplicity, we assume that $\mathcal{I} = \{1, 2, \ldots, \mathbf{I}\}$ and then $|\mathcal{I}| = \mathbf{I}$. For each $i \in \mathcal{I}$, let $\gamma^i(\cdot)$ be a vector valued function defined from $\mathbb{R}^d$ into $\mathbb{R}^d$. Fix $\mu \in \mathbb{R}^d$, $\Gamma$ a $d \times d$ symmetric and strictly positive definite covariance matrix and $\nu$ a probability measure on $(\overline{G}, \mathcal{B}(\overline{G}))$, where $\mathcal{B}(\overline{G})$ denotes the $\sigma$-algebra of Borel subsets of the closure $\overline{G}$ of $G$.

DEFINITION 2.1 (*Semimartingale reflecting Brownian motion*). A semimartingale reflecting Brownian motion (abbreviated as SRBM) associated with the data $(G, \mu, \Gamma, \{\gamma^i, i \in \mathcal{I}\}, \nu)$ is an $\{\mathcal{F}_t\}$-adapted, $d$-dimensional process $W$ defined on some filtered probability space $(\Omega, \mathcal{F}, \{\mathcal{F}_t\}, P)$ such that:

(i) $P$-a.s., $W(t) = X(t) + \sum_{i \in \mathcal{I}} \int_{(0,t]} \gamma^i(W(s)) \, dY_i(s)$ for all $t \geq 0$,

(ii) $P$-a.s., $W$ has continuous paths and $W(t) \in \overline{G}$ for all $t \geq 0$,

(iii) under $P$, $X$ is a $d$-dimensional $\{\mathcal{F}_t\}$-Brownian motion with drift vector $\mu$, covariance matrix $\Gamma$ and initial distribution $\nu$,

(iv) for each $i \in \mathcal{I}$, $Y_i$ is an $\{\mathcal{F}_t\}$-adapted, one-dimensional process such that $P$-a.s.,

(a) $Y_i(0) = 0$,

(b) $Y_i$ is continuous and nondecreasing,

(c) $Y_i(t) = \int_{(0,t]} 1_{\{W(s) \in \partial G_i \cap \partial G\}} dY_i(s)$ for all $t \geq 0$.

We shall often refer to $Y = \{Y_i, i \in \mathcal{I}\}$ as the "pushing process" associated with the SRBM $W$. When $\nu = \delta_x$, we may alternatively say that $W$ is an SRBM associated with the data $(G, \mu, \Gamma, \{\gamma^i, i \in \mathcal{I}\})$ that starts from $x$. We will call $(W, X, Y)$ satisfying Definition 2.1 an *extended* SRBM associated with the data $(G, \mu, \Gamma, \{\gamma^i, i \in \mathcal{I}\}, \nu)$.

Loosely speaking, an SRBM behaves like a Brownian motion in the interior of the domain $G$ and it is confined to $\overline{G}$ by instantaneous "reflection" (or "pushing") at the boundary, where the allowed directions of "reflection" at $x \in \partial G$ are convex combinations of the vectors $\gamma^i(x)$ for $i$ such that $x \in \partial G_i$. Under the assumptions imposed on $G$ and $\{\gamma^i, i \in \mathcal{I}\}$ in Sections 3.1 and 3.2 below, at each point on the boundary of $G$ there is an allowed direction of reflection that can be used there which "points into the interior of $G$." We end this section by introducing a related set-valued function $\mathcal{I}(\cdot)$ and show a key property of it.

DEFINITION 2.2. For each $x \in \mathbb{R}^d$, let $\mathcal{I}(x) = \{i \in \mathcal{I} : x \in \partial G_i\}$.

The set-valued function $\mathcal{I}(\cdot)$ has the following property called upper semicontinuity on $\partial G$.

LEMMA 2.1. *For each $x \in \partial G$, there is an open neighborhood $V_x$ of $x$ in $\mathbb{R}^d$ such that*

(4) $$\mathcal{I}(y) \subset \mathcal{I}(x) \quad \text{for all } y \in V_x.$$

PROOF. We prove this lemma by contradiction. Suppose that the function $\mathcal{I}(\cdot)$ does not satisfy (4). Then there is a point $x \in \partial G$ such that there is no open neighborhood $V_x$ of $x$ such that $\mathcal{I}(y) \subset \mathcal{I}(x)$ for all $y \in V_x$. Since the index set $\mathcal{I}$ is finite, there is an index $k \in \mathcal{I} \setminus \mathcal{I}(x)$ and a sequence of points $\{y_n\} \subset \mathbb{R}^d$ such that $\|y_n - x\| < \frac{1}{n}$ and $k \in \mathcal{I}(y_n)$ for each $n \geq 1$. Hence $y_n \in \partial G_k$ for all $n \geq 1$. Since $\partial G_k$ is closed and $y_n \to x$ as $n \to \infty$, we conclude that $x \in \partial G_k$. This implies that $k \in \mathcal{I}(x)$, which is a contradiction, as desired. $\square$

## 3. Assumptions on the domain $G$ and the reflection vector fields $\{\gamma^i\}$.

3.1. *Assumptions on the domain $G$.* We henceforth assume that the domain $G$ satisfies assumptions (A1)–(A3) below. In the case when $G$ is bounded, assumptions (A2)–(A3) follow from assumption (A1) (see Lemmas A.1 and A.2 in the Appendix for details). If the domain $G$ is a convex polyhedron satisfying assumption (A1), then assumptions (A2)–(A3) hold by Lemma A.3 in the Appendix.

(A1) $G$ is a nonempty domain in $\mathbb{R}^d$ with representation

$$G = \bigcap_{i \in \mathcal{I}} G_i, \tag{5}$$

where for each $i \in \mathcal{I}$, $G_i$ is a nonempty domain, $G_i \neq \mathbb{R}^d$, and the boundary $\partial G_i$ of $G_i$ is $C^1$. For each $i \in \mathcal{I}$, we let $n^i(\cdot)$ be the unit normal vector field on $\partial G_i$ that points into $G_i$.

(A2) For each $\varepsilon \in (0,1)$ there exists $R(\varepsilon) > 0$ such that for each $i \in \mathcal{I}$, $x \in \partial G_i \cap \partial G$ and $y \in \overline{G}$ satisfying $\|x - y\| < R(\varepsilon)$, we have

$$\langle n^i(x), y - x \rangle \geq -\varepsilon \|x - y\|. \tag{6}$$

(A3) The function $D : [0, \infty) \to [0, \infty]$ defined such that $D(0) = 0$ and

$$D(r) = \sup_{\substack{\mathcal{J} \subset \mathcal{I} \\ \mathcal{J} \neq \varnothing}} \sup \left\{ \operatorname{dist}\left(x, \bigcap_{j \in \mathcal{J}} (\partial G_j \cap \partial G)\right) : x \in \bigcap_{j \in \mathcal{J}} U_r(\partial G_j \cap \partial G) \right\} \tag{7}$$

for $r > 0$, satisfies

$$D(r) \to 0 \quad \text{as } r \to 0. \tag{8}$$

REMARK. Assumption (A2) is reminiscent of the uniform exterior cone condition (cf. [9], page 195). We say that a region $G \subset \mathbb{R}^d$ satisfies a uniform exterior cone condition if for each $x_0 \in \partial G$, there is a truncated closed right circular cone $V_{x_0}$, with nonempty interior and vertex $x_0$, satisfying $V_{x_0} \cap \overline{G} = \{x_0\}$, and the truncated cones $V_{x_0}$ are all congruent to some fixed truncated closed right circular cone $V$. By comparing assumption (A2) with the uniform exterior cone condition, we see that assumption (A2) implies the uniform exterior cone condition. On the other hand, under assumption (A1), assumption (A2) is implied by a family of uniform exterior cone conditions where for each $\varepsilon \in (0,1)$, the axis of the truncated closed right circular cone at $x \in \partial G$ is along the vector $-n^i(x)$ and all of the truncated closed right circular cones are congruent to a truncated closed right circular cone whose height and base radius are $R(\varepsilon)$ and $R(\varepsilon)(\frac{1}{\varepsilon^2} - 1)^{1/2}$ respectively. Assumption (A2) holds automatically if $G$ is convex. We also note that assumption (A2) is strictly weaker than the uniform exterior sphere condition. The definition of the uniform exterior sphere condition is similar to that for the uniform exterior cone condition where a closed ball with $x_0$ on its boundary takes the place of the truncated closed right circular cone $V_{x_0}$. It can be checked that for the domain $G = \{(x,y) \in \mathbb{R}^2 : y < |x|^\alpha\}$ with $\alpha \in (1,2)$, the uniform exterior sphere condition fails to hold, but assumption (A2) holds. In fact, at the point $(0,0) \in \mathbb{R}^2$, there is no $r > 0$ and $y \in \mathbb{R}^2$ such that $B_r(y) \cap \partial G = \{(0,0)\}$.

REMARK. For the definition of $D(\cdot)$ in (A3), we adopt the convention that the supremum over an empty set is zero and $\operatorname{dist}(x,\varnothing)=\infty$. Since $\partial G_i \cap \partial G \neq \varnothing$ for at least one $i \in \mathcal{I}$, the function $D(\cdot)$ satisfies $\lim_{r\to\infty} D(r) = \infty$. Furthermore, $D(r_1) \leq D(r_2)$ whenever $r_1, r_2 \in [0,\infty)$ and $r_1 \leq r_2$. Assumption (A3) requires that for any nonempty subset $\mathcal{J} \subset \mathcal{I}$, the intersection of tubular neighborhoods of the boundaries $\{\partial G_j \cap \partial G : j \in \mathcal{J}\}$ given by the set $\bigcap_{j \in \mathcal{J}} U_r(\partial G_j \cap \partial G)$ "converges" to the intersection of the boundaries given by the set $\bigcap_{j \in \mathcal{J}} (\partial G_j \cap \partial G)$ as $r$ approaches 0. Property (8) need not always hold. For example, let $G_1 = \{(x,y) \in \mathbb{R}^2 : y < e^{-x^2/2}, x \in \mathbb{R}\}$ and $G_2 = \{(x,y) \in \mathbb{R}^2 : y > 0, x \in \mathbb{R}\}$. Then $\partial G_1 \cap \partial G_2 = \varnothing$. But for each $r > 0$, $U_r(\partial G_1) \cap U_r(\partial G_2) \neq \varnothing$. Hence $D(r) = \infty$ for each $r > 0$.

3.2. *Assumptions on the reflection vector fields $\{\gamma^i\}$.* We henceforth assume that there are vector fields $\{\gamma^i(\cdot), i \in \mathcal{I}\}$ satisfying assumptions (A4)–(A5) below.

(A4) There is a constant $L > 0$ such that for each $i \in \mathcal{I}$, $\gamma^i(\cdot)$ is a uniformly Lipschitz continuous function from $\mathbb{R}^d$ into $\mathbb{R}^d$ with Lipschitz constant $L$ and $\|\gamma^i(x)\| = 1$ for each $x \in \mathbb{R}^d$.

(A5) There is a constant $a \in (0,1)$, and vector valued functions $b(\cdot) = (b_1(\cdot), \ldots, b_\mathbf{I}(\cdot))$ and $c(\cdot) = (c_1(\cdot), \ldots, c_\mathbf{I}(\cdot))$ from $\partial G$ into $\mathbb{R}_+^\mathbf{I}$ such that for each $x \in \partial G$,

(i) $\sum_{i \in \mathcal{I}(x)} b_i(x) = 1$,

$$(9) \qquad \min_{j \in \mathcal{I}(x)} \left\langle \sum_{i \in \mathcal{I}(x)} b_i(x) n^i(x), \gamma^j(x) \right\rangle \geq a,$$

(ii) $\sum_{i \in \mathcal{I}(x)} c_i(x) = 1$,

$$(10) \qquad \min_{j \in \mathcal{I}(x)} \left\langle \sum_{i \in \mathcal{I}(x)} c_i(x) \gamma^i(x), n^j(x) \right\rangle \geq a.$$

We note here for future use that by (A4), if we set $\rho_0 = \frac{a}{4L}$, then for any $x, y \in \mathbb{R}^d$ satisfying $\|x - y\| < \rho_0$, we have $\|\gamma^i(x) - \gamma^i(y)\| < a/4$ for each $i \in \mathcal{I}$. So for each $0 < \rho < \rho_0/4$, by (9)–(10) and the normalization of $b(\cdot), c(\cdot), \gamma^i(\cdot), n^j(\cdot)$ for $i, j \in \mathcal{I}$, we obtain

$$(11) \qquad \inf_{x \in \partial G} \min_{j \in \mathcal{I}(x)} \inf_{y \in B_{4\rho}(x)} \left\langle \sum_{i \in \mathcal{I}(x)} b_i(x) n^i(x), \gamma^j(y) \right\rangle \geq a/2$$

and

$$(12) \qquad \inf_{x \in \partial G} \min_{j \in \mathcal{I}(x)} \inf_{y \in B_{4\rho}(x)} \left\langle \sum_{i \in \mathcal{I}(x)} c_i(x) \gamma^i(y), n^j(x) \right\rangle \geq a/2.$$

The use of $B_{4\rho}(x)$ here is related to the form in which this is used in Section 4.1.

REMARK. Assumption (A4) is equivalent to (3.4) in [6] when $G$ is bounded. Property (10) means that, at each point $x \in \partial G$, there is a convex combination $\gamma(x) = \sum_{i \in \mathcal{I}(x)} c_i(x) \gamma^i(x)$ of the vectors $\{\gamma^i(x), i \in \mathcal{I}(x)\}$ that can be used there such that $\gamma(x)$ "points into" $G$. Property (9) is in a sense a dual condition to property (10), where the roles of $\gamma^i$ and $n^i$ are reversed for $i \in \mathcal{I}(x)$. This property (9) is used in showing the oscillation inequality in Theorem 4.1 below. Assumption (A5) is an analogue of Assumption 1.1 in [4]. When $G$ is bounded, (10) is similar to condition (3.6) in [6] (we assume some additional uniformity through the lack of dependence of $a$ on $x$).

It is straightforward to see using the triangle inequality that the following condition (A5)$'$ implies (A5).

(A5)$'$ There is $a \in (0,1)$ and vector valued functions $b, c$ from $\partial G$ into $\mathbb{R}_+^{\mathbf{I}}$ such that for each $x \in \partial G$,

(i) $\sum_{i \in \mathcal{I}(x)} b_i(x) = 1$, and for each $i \in \mathcal{I}(x)$,

$$(13) \qquad b_i(x) \langle n^i(x), \gamma^i(x) \rangle \geq a + \sum_{j \in \mathcal{I}(x) \setminus \{i\}} b_j(x) |\langle n^j(x), \gamma^i(x) \rangle|,$$

(ii) $\sum_{i \in \mathcal{I}(x)} c_i(x) = 1$, and for each $i \in \mathcal{I}(x)$,

$$(14) \qquad c_i(x) \langle \gamma^i(x), n^i(x) \rangle \geq a + \sum_{j \in \mathcal{I}(x) \setminus \{i\}} c_j(x) |\langle \gamma^j(x), n^i(x) \rangle|.$$

Condition (A5)$'$(ii) is similar to condition (3.8) in [6], although here we assume additional uniformity through the lack of dependence of $a$ on $x$. As noted in [6], their condition (3.8) can be phrased in terms of a nonsingular M-matrix requirement [2]. (This is sometimes also called a generalized Harrison–Reiman type of condition [10].) Since that nonsingular M-matrix property is invariant under transpose, and this property for the transpose corresponds to a local form of (A5)$'$(i), one might conjecture that there is an equivalence between the existence of a nonnegative vector valued function $b$ such that (A5)$'$(i) holds for each $x \in \partial G$ and the existence of a nonnegative vector valued function $c$ such that (A5)$'$(ii) holds for each $x \in \partial G$. Indeed we have the following lemma. We have stated the two (equivalent) conditions (i) and (ii) in specifying (A5)$'$ to preserve a parallel with (A5) and since both properties can be useful in proofs. Furthermore, in light of the following lemma, verifying either condition suffices for both to hold.

LEMMA 3.1. *There is a constant $a \in (0,1)$ and a vector valued function $b : \partial G \to \mathbb{R}_+^{\mathbf{I}}$ such that $(A5)'(i)$ holds for each $x \in \partial G$ if and only if there is a constant $a \in (0,1)$ and a vector valued function $c : \partial G \to \mathbb{R}_+^{\mathbf{I}}$ such that $(A5)'(ii)$ holds for each $x \in \partial G$.*

PROOF. We just prove the "if" part; the "only if" part can be proved in a similar manner.

We suppose that there is a constant $a \in (0,1)$ and a vector valued function $c : \partial G \to \mathbb{R}_+^{\mathbf{I}}$ such that $(A5)'(ii)$ holds for each $x \in \partial G$. For fixed $x \in \partial G$, consider the square matrix $A(x)$ whose diagonal entries are given by the (positive) elements $\langle n^i(x), \gamma^i(x) \rangle$ for $i \in \mathcal{I}(x)$ and whose off-diagonal entries are given by $-|\langle n^i(x), \gamma^j(x) \rangle|$ for $i \in \mathcal{I}(x), j \in \mathcal{I}(x), j \neq i$. Let $E$ be the square matrix having the same dimensions as $A(x)$ and whose entries are all equal to one. By the theory of M-matrices (see [2], Chapter 6, especially condition $(M_{35})$), condition (ii) of $(A5)'$ implies that $A(x) - \frac{a}{2}E$ is a nonsingular M-matrix, that is, $A(x) - \frac{a}{2}E$ has nonnegative diagonal entries and nonpositive off-diagonal entries and it can be written in the form $s(x)I - B(x)$ where $B(x)$ is a matrix with nonnegative entries and $s(x) > 0$ is a constant that is strictly larger than the spectral radius of $B(x)$.

Since the nonsingular M-matrix property is invariant under transpose (cf. $(G_{21})$ in Chapter 6 of [2]), then $A'(x) - \frac{a}{2}E$ is also a nonsingular M-matrix. Hence, there is a vector $\tilde{b}(x) = (\tilde{b}_i(x) : i \in \mathcal{I}(x))$ with nonnegative entries such that $(A'(x) - \frac{a}{2}E)\tilde{b}(x) > 0$ (cf. $(I_{27})$ in Chapter 6 of [2]). We can extend $\tilde{b}(x)$ to an **I**-dimensional vector $b(x)$ and normalize it so that $\sum_{i \in \mathcal{I}(x)} b_i(x) = 1$. Then $(A5)'(i)$ holds with $\frac{a}{2}$ in place of $a$. □

**4. Invariance principle.** In this section we state and prove an invariance principle for an SRBM living in the closure of a domain $G$ with piecewise smooth boundary and having associated reflection fields $\{\gamma^i, i \in \mathcal{I}\}$, where $G$, $\{\gamma^i, i \in \mathcal{I}\}$ satisfy assumptions (A1)–(A5) of Section 3. (These assumptions hold throughout this section.) We shall first state a preliminary result called an *oscillation inequality* (see Theorem 4.1), then we use it to prove a tightness result (see Theorem 4.2). Finally, we establish the invariance principle (see Theorem 4.3).

4.1. *Oscillation inequality.* The following oscillation inequality is the key to the proof of the tightness result claimed in Theorem 4.2. In this subsection, for any $0 \leq t_1 < t_2 < \infty$ and any integer $k \geq 1$, $D([t_1, t_2], \mathbb{R}^k)$ denotes the set of functions $w : [t_1, t_2] \to \mathbb{R}^k$ that are right continuous on $[t_1, t_2)$ and have finite left limits on $(t_1, t_2]$. For $w \in D([t_1, t_2], \mathbb{R}^k)$,

(15) $\quad \mathrm{Osc}(w, [t_1, t_2]) = \sup\{\|w(t) - w(s)\| : t_1 \leq s < t \leq t_2\},$

(16) $\quad \mathrm{Osc}(w, [t_1, t_2)) = \sup\{\|w(t) - w(s)\| : t_1 \leq s < t < t_2\}.$

Note that we do not explicitly indicate the dependence on $k$ in the notation.

Recall the constants $a, L$ from assumptions (A4)–(A5), the functions $R(\cdot)$ from assumption (A2) and $D(\cdot)$ from (7). Let $\rho_0 = \frac{a}{4L}$.

THEOREM 4.1 (Oscillation inequality). *There exists a nondecreasing function $\Pi : (0, \infty) \to (0, \infty]$ satisfying $\Pi(u) \to 0$ as $u \to 0$, such that $\Pi$ depends only on the constants $\mathbf{I}, a$ and the function $D(\cdot)$, and such that whenever $0 < \rho < \min\{\frac{\rho_0}{4}, \frac{R(a/4)}{4}\}$, $0 < \delta < \frac{\rho}{2}$, $0 \le s < t < \infty$, $w, x \in D([s,t], \mathbb{R}^d)$ and $y \in D([s,t], \mathbb{R}^{\mathbf{I}})$ satisfy:*

(i) $w(u) \in B_\rho(x_0) \cap U_\delta(G)$ *for all* $u \in [s,t]$, *for some* $x_0 \in \overline{G}$,
(ii) $w(u) = w(s) + x(u) - x(s) + \sum_{i \in \mathcal{I}} \int_{(s,u]} \gamma^i(w(v)) \, dy_i(v)$ *for all* $u \in [s,t]$,
(iii) *for each* $i \in \mathcal{I}$,

    (a) $y_i(s) \ge 0$,
    (b) $y_i$ *is nondecreasing and* $\Delta y_i(u) \le \delta$ *for all* $u \in (s,t]$,
    (c) $y_i(u) = y_i(s) + \int_{(s,u]} 1_{\{w(v) \in U_\delta(\partial G_i \cap \partial G)\}} \, dy_i(v)$ *for all* $u \in [s,t]$,

(iv) $D(\Pi(\operatorname{Osc}(x, [s,t]) + \delta)) < \frac{\rho}{2}$,

*then we have that the following hold:*

(17) $$\operatorname{Osc}(w, [s,t]) \le \Pi(\operatorname{Osc}(x, [s,t]) + \delta),$$

(18) $$\operatorname{Osc}(y, [s,t]) \le \Pi(\operatorname{Osc}(x, [s,t]) + \delta).$$

PROOF. Let

$$\Pi_0(u) = u \qquad \text{for all } u > 0.$$

Define $\Pi_m : (0, \infty) \to (0, \infty]$, $m = 1, \ldots, \mathbf{I}$, inductively such that

$$\Pi_m(u) = \Pi_{m-1}(u) + (\mathbf{I} + 2)u + \left(1 + \frac{4}{a}\right)(D(\Pi_{m-1}(u) + (\mathbf{I} + 2)u) + 2u).$$

Here the sum of any element of $[0, \infty)$ with $\infty$ is $\infty$ and $D(\infty)$ is defined to equal $\infty$. For each $m = 0, 1, \ldots, \mathbf{I}$, the function $\Pi_m$ is nondecreasing and depends only on $\mathbf{I}, a$ and $D(\cdot)$. For each $m = 1, \ldots, \mathbf{I}$ and $u > 0$, $\Pi_{m-1}(u) \le \Pi_m(u)$. By assumption (A3), we conclude (using an induction proof) that

$$\Pi_m(u) \to 0 \quad \text{as } u \to 0, \text{ for } m = 0, 1, \ldots, \mathbf{I}.$$

Let $\Pi(\cdot) = \Pi_{\mathbf{I}}(\cdot)$.

Fix $0 < \rho < \min\{\frac{\rho_0}{4}, \frac{R(a/4)}{4}\}$, $0 < \delta < \frac{\rho}{2}$, $0 \le s < t < \infty$. Suppose that $w, x \in D([s,t], \mathbb{R}^d)$ and $y \in D([s,t], \mathbb{R}^{\mathbf{I}})$ satisfy (i)–(iv) in the statement of Theorem 4.1. For each nonempty interval $[t_1, t_2] \subset [s,t]$, let

$$\mathcal{I}_{[t_1, t_2]} = \{i \in \mathcal{I} : w(u) \in U_\delta(\partial G_i \cap \partial G) \text{ for some } u \in [t_1, t_2]\},$$

the indices of the boundary surfaces that $w(\cdot)$ comes close to in the time interval $[t_1, t_2]$. For each $0 \leq m \leq \mathbf{I}$, define $\mathcal{T}_m = \{[t_1, t_2] \subset [s, t] : |\mathcal{I}_{[t_1, t_2]}| \leq m\}$. Note that under the partial ordering of set inclusion, $\mathcal{T}_m$ increases with $m$. To prove the theorem, we will prove by induction that for each $0 \leq m \leq \mathbf{I}$ and each interval $[t_1, t_2] \in \mathcal{T}_m$, (17)–(18) hold with $[t_1, t_2]$ in place of $[s, t]$ and $\Pi_m(\cdot)$ in place of $\Pi(\cdot)$. The result for $m = \mathbf{I}$ yields the theorem.

Suppose that $m = 0$. Then $\mathcal{T}_0 = \{[t_1, t_2] \subset [s, t] : |\mathcal{I}_{[t_1, t_2]}| = 0\}$. Fix an interval $[t_1, t_2] \in \mathcal{T}_0$. Since $\mathcal{I}_{[t_1, t_2]} = \varnothing$ and (iii)(c) holds, the function $y$ does not increase on the time interval $(t_1, t_2]$, that is, $y_i(t_2) - y_i(t_1) = 0$ for all $i \in \mathcal{I}$. Then, for $t_1 \leq u < v \leq t_2$,

$$\text{(19)} \qquad w(v) - w(u) = x(v) - x(u).$$

So in this case,

$$\text{(20)} \qquad \operatorname{Osc}(w, [t_1, t_2]) = \operatorname{Osc}(x, [t_1, t_2]) \leq \operatorname{Osc}(x, [t_1, t_2]) + \delta,$$

$$\text{(21)} \qquad \operatorname{Osc}(y, [t_1, t_2]) = 0 \leq \operatorname{Osc}(x, [t_1, t_2]) + \delta.$$

Thus, (17)–(18) hold with $\Pi_0(\cdot)$ in place of $\Pi(\cdot)$ and $[t_1, t_2]$ in place of $[s, t]$ for each interval $[t_1, t_2] \in \mathcal{T}_0$.

For the induction step, let $1 \leq m \leq \mathbf{I}$ and suppose that (17)–(18) hold with $\Pi_{m-1}(\cdot)$ in place of $\Pi(\cdot)$ and $[t_1, t_2]$ in place of $[s, t]$ for each interval $[t_1, t_2] \in \mathcal{T}_{m-1}$.

Now fix $[t_1, t_2] \in \mathcal{T}_m$. If $|\mathcal{I}_{[t_1, t_2]}| \leq m - 1$, then $[t_1, t_2] \in \mathcal{T}_{m-1}$ and so by the induction assumption we have that (17)–(18) hold with $[t_1, t_2]$ in place of $[s, t]$ and $\Pi_{m-1}(\cdot)$ [and hence $\Pi_m(\cdot)$] in place of $\Pi(\cdot)$. Thus, it suffices to consider $[t_1, t_2] \subset [s, t]$ such that $|\mathcal{I}_{[t_1, t_2]}| = m$. For $i \notin \mathcal{I}_{[t_1, t_2]}$, by (iii)(c), $y_i(t_2) - y_i(t_1) = 0$, and so by (ii), for $t_1 \leq u < v \leq t_2$, we have

$$\text{(22)} \qquad w(v) - w(u) = x(v) - x(u) + \sum_{i \in \mathcal{I}_{[t_1, t_2]}} \int_{(u, v]} \gamma^i(w(r)) \, dy_i(r).$$

Let $\overline{\Pi}_m(u) = \Pi_{m-1}(u) + (\mathbf{I} + 2)u$ for all $u > 0$, and $\eta = \operatorname{Osc}(x, [t_1, t_2]) + \delta$. For any $M \in (0, \infty]$ and any nonempty set $\mathcal{J} \subset \mathcal{I}$, let

$$F_{\mathcal{J}}^M = \{z \in \mathbb{R}^d : \operatorname{dist}(z, \partial G_i \cap \partial G) < M \text{ for all } i \in \mathcal{J}\}.$$

Note that $F_{\mathcal{J}}^M = \varnothing$ when there is an $i \in \mathcal{J}$ such that $\partial G_i \cap \partial G = \varnothing$. Since $\overline{\Pi}_m(\cdot) \leq \Pi_m(\cdot) \leq \Pi(\cdot)$, $D(\cdot)$ and $\Pi(\cdot)$ are nondecreasing, and $\operatorname{Osc}(x, [t_1, t_2]) \leq \operatorname{Osc}(x, [s, t])$, we have by (iv) that

$$\text{(23)} \qquad D(\overline{\Pi}_m(\eta)) \leq D(\Pi_m(\eta)) \leq D(\Pi(\eta)) < \frac{\rho}{2}.$$

Note that this implies $\overline{\Pi}_m(\eta) < \infty$ since $D(\infty) = \infty$.

We now consider two cases.

*Case* 1. Suppose that $w(r) \in F^{\overline{\Pi}_m(\eta)}_{\mathcal{I}_{[t_1,t_2]}}$ for all $r \in [t_1, t_2]$.

Fix $u, v$ such that $t_1 \leq u < v \leq t_2$. Since we have that
$$w(v) \in \bigcap_{j \in \mathcal{I}_{[t_1,t_2]}} U_{\overline{\Pi}_m(\eta)}(\partial G_j \cap \partial G),$$
by the definition of $D(\cdot)$ and (23), there is $z \in \bigcap_{j \in \mathcal{I}_{[t_1,t_2]}}(\partial G_j \cap \partial G)$ such that
$$\text{(24)} \qquad \|w(v) - z\| \leq D(\overline{\Pi}_m(\eta)) < \frac{\rho}{2}.$$

For each $r \in [t_1, t_2]$, by (i) we have that $w(r) \in U_\delta(G)$, and so there is $z^r \in G$ such that
$$\|w(r) - z^r\| \leq 2\delta.$$

Hence by (i) and (24) we have
$$\text{(25)} \qquad \begin{aligned} \|z^r - z\| &\leq \|z^r - w(r)\| + \|w(r) - x_0\| + \|x_0 - w(v)\| + \|w(v) - z\| \\ &\leq 2\delta + \rho + \rho + \rho/2 < 4\rho < R(a/4) \end{aligned}$$

and
$$\text{(26)} \qquad \begin{aligned} \|w(r) - z\| &\leq \|w(r) - x_0\| + \|x_0 - w(v)\| + \|w(v) - z\| \\ &\leq \rho + \rho + \rho/2 < 4\rho. \end{aligned}$$

By (6) and (25) we have
$$\text{(27)} \quad \langle n^j(z), z - z^r \rangle \leq \frac{a}{4}\|z - z^r\| \qquad \text{for each } j \in \mathcal{I}(z) \text{ and } r \in [t_1, t_2].$$

Note that $\mathcal{I}(z) \supset \mathcal{I}_{[t_1,t_2]}$. Recalling the definition of $b(\cdot)$ from assumption (A5), on dotting the vector $\sum_{j \in \mathcal{I}(z)} b_j(z) n^j(z)$ with both sides of (22) and rearranging, we obtain
$$\text{(28)} \qquad \begin{aligned} \sum_{i \in \mathcal{I}_{[t_1,t_2]}} \int_{(u,v]} &\left\langle \sum_{j \in \mathcal{I}(z)} b_j(z) n^j(z), \gamma^i(w(r)) \right\rangle dy_i(r) \\ &= \sum_{j \in \mathcal{I}(z)} b_j(z) \langle n^j(z), w(v) - w(u) \rangle \\ &\quad - \sum_{j \in \mathcal{I}(z)} b_j(z) \langle n^j(z), x(v) - x(u) \rangle. \end{aligned}$$

So by (11), (22), (24)–(28), and the fact that $\sum_{j \in \mathcal{I}(z)} b_j(z) = 1$, $b_j(z) \geq 0$ for $j \in \mathcal{I}$, we have
$$\frac{a}{2} \sum_{i \in \mathcal{I}_{[t_1,t_2]}} (y_i(v) - y_i(u))$$

$$\leq \sum_{i\in\mathcal{I}_{[t_1,t_2]}} \int_{(u,v]} \left\langle \sum_{j\in\mathcal{I}(z)} b_j(z) n^j(z), \gamma^i(w(r)) \right\rangle dy_i(r)$$

$$= \sum_{j\in\mathcal{I}(z)} b_j(z)\langle n^j(z), w(v) - z\rangle + \sum_{j\in\mathcal{I}(z)} b_j(z)\langle n^j(z), z - z^u\rangle$$

$$+ \sum_{j\in\mathcal{I}(z)} b_j(z)\langle n^j(z), z^u - w(u)\rangle - \sum_{j\in\mathcal{I}(z)} b_j(z)\langle n^j(z), x(v) - x(u)\rangle$$

$$\leq D(\overline{\Pi}_m(\eta)) + \frac{a}{4}\|z - z^u\| + 2\delta + \|x(v) - x(u)\|$$

$$\leq D(\overline{\Pi}_m(\eta)) + 2\delta + \|x(v) - x(u)\|$$
$$+ \frac{a}{4}(\|z - w(v)\| + \|w(v) - w(u)\| + \|w(u) - z^u\|)$$

$$\leq D(\overline{\Pi}_m(\eta)) + 2\delta + \|x(v) - x(u)\|$$
$$+ \frac{a}{4}\left(D(\overline{\Pi}_m(\eta)) + \|x(v) - x(u)\| + \sum_{i\in\mathcal{I}_{[t_1,t_2]}} (y_i(v) - y_i(u)) + 2\delta\right)$$

$$\leq \left(1 + \frac{a}{4}\right)\{D(\overline{\Pi}_m(\eta)) + 2\delta + \|x(v) - x(u)\|\} + \frac{a}{4}\sum_{i\in\mathcal{I}_{[t_1,t_2]}} (y_i(v) - y_i(u)).$$

Hence

$$\frac{a}{4}\sum_{i\in\mathcal{I}_{[t_1,t_2]}} (y_i(v) - y_i(u)) \leq \left(1 + \frac{a}{4}\right)\{D(\overline{\Pi}_m(\eta)) + 2\delta + \|x(v) - x(u)\|\}$$

$$\leq \left(1 + \frac{a}{4}\right)\{D(\overline{\Pi}_m(\eta)) + 2\eta\}.$$

On multiplying through by $\frac{4}{a}$, we obtain

$$(29) \quad \sum_{i\in\mathcal{I}_{[t_1,t_2]}} (y_i(v) - y_i(u)) \leq \left(1 + \frac{4}{a}\right)\{D(\overline{\Pi}_m(\eta)) + 2\eta\} \leq \Pi_m(\eta).$$

Hence, by (29) and the fact that for any $x \in \mathbb{R}^d$, $\|x\| \leq \sum_{i=1}^d |x_i|$, we have

$$(30) \quad \mathrm{Osc}(y, [t_1, t_2]) \leq \Pi_m(\mathrm{Osc}(x, [t_1, t_2]) + \delta),$$

and by (22), (29) and the definitions of $\overline{\Pi}_m(\cdot)$ and $\Pi_m(\cdot)$, we have

$$\mathrm{Osc}(w, [t_1, t_2]) \leq \mathrm{Osc}(x, [t_1, t_2]) + \left(1 + \frac{4}{a}\right)\{D(\overline{\Pi}_m(\eta)) + 2\eta\}$$
$$\leq \Pi_m(\mathrm{Osc}(x, [t_1, t_2]) + \delta),$$

as desired.

*Case* 2. Suppose that there is $t_3 \in [t_1, t_2]$ such that $w(t_3) \notin F_{\mathcal{I}_{[t_1,t_2]}}^{\overline{\Pi}_m(\eta)}$.

Define $\sigma = \inf\{u \in [t_1, t_2] : w(u) \notin F_{\mathcal{I}_{[t_1,t_2]}}^{\overline{\Pi}_m(\eta)}\}$. Then $\sigma \leq t_2$. For each $u \in [t_1, \sigma)$, $w(u) \in F_{\mathcal{I}_{[t_1,t_2]}}^{\overline{\Pi}_m(\eta)}$ and so by a similar analysis to that for Case 1, we obtain for each $v \in [t_1, \sigma)$,

$$\mathrm{Osc}(w, [t_1, v]) \leq \eta + \left(1 + \frac{4}{a}\right)(D(\overline{\Pi}_m(\eta)) + 2\eta)$$

and

$$\mathrm{Osc}(y, [t_1, v]) \leq \left(1 + \frac{4}{a}\right)(D(\overline{\Pi}_m(\eta)) + 2\eta).$$

By the right continuity of paths we have $w(\sigma) \notin F_{\mathcal{I}_{[t_1,t_2]}}^{\overline{\Pi}_m(\eta)}$. Then there is an $i \in \mathcal{I}_{[t_1,t_2]}$ such that $\mathrm{dist}(w(\sigma), \partial G_i \cap \partial G) \geq \overline{\Pi}_m(\eta)$, and it follows that $w$ does not reach $U_\delta(\partial G_i \cap \partial G)$ during the interval $[\sigma, t_2]$. To see this, let $\tau = \inf\{u \in [\sigma, t_2] : \mathrm{dist}(w(u), \partial G_i \cap \partial G) \leq \delta\}$ with the convention that the infimum of an empty set is $\infty$. If $\tau \leq t_2$, then by the right continuity of $w(\cdot)$ and since $\overline{\Pi}_m(\eta) > \delta$, we have $\tau > \sigma$ and $\mathrm{dist}(w(\tau), \partial G_i \cap \partial G) \leq \delta$. Also, since $|\mathcal{I}_{[t_1,t_2]}| = m$, we have $[\sigma, u] \in \mathcal{T}_{m-1}$ for each $u \in [\sigma, \tau)$. By the induction assumption and letting $u \to \tau$, we have $\|w(\tau-) - w(\sigma)\| \leq \Pi_{m-1}(\eta)$. By (ii), (iii)(b) and since $\|\gamma^i(\cdot)\| = 1$, we have

$$\|\Delta w(\tau)\| \leq \|\Delta x(\tau)\| + \sum_{i \in \mathcal{I}} \Delta y_i(\tau) \leq \mathrm{Osc}(x, [t_1, t_2]) + \mathbf{I}\delta \leq \mathbf{I}\eta.$$

Then simple manipulations yield

$$\mathrm{dist}(w(\sigma), \partial G_i \cap \partial G) \leq \|w(\sigma) - w(\tau-)\| + \|\Delta w(\tau)\| + \mathrm{dist}(w(\tau), \partial G_i \cap \partial G)$$
$$\leq \Pi_{m-1}(\eta) + \mathbf{I}\eta + \delta$$
$$< \overline{\Pi}_m(\eta).$$

This contradicts the fact that $\mathrm{dist}(w(\sigma), \partial G_i \cap \partial G) \geq \overline{\Pi}_m(\eta)$, and so confirms that $w$ does not reach $U_\delta(\partial G_i \cap \partial G)$ in $[\sigma, t_2]$. Thus we must have $[\sigma, t_2] \in \mathcal{T}_{m-1}$. Hence we have by the induction assumption that

$$\mathrm{Osc}(w, [t_1, t_2]) \leq \sup_{v \in [t_1, \sigma)} \mathrm{Osc}(w, [t_1, v]) + \|\Delta w(\sigma)\| + \mathrm{Osc}(w, [\sigma, t_2])$$
$$\leq \eta + \left(1 + \frac{4}{a}\right)(D(\overline{\Pi}_m(\eta)) + 2\eta) + \mathbf{I}\eta + \Pi_{m-1}(\eta)$$
$$\leq \Pi_m(\mathrm{Osc}(x, [t_1, t_2]) + \delta)$$

and

$$\mathrm{Osc}(y, [t_1, t_2]) \leq \sup_{v \in [t_1, \sigma)} \mathrm{Osc}(y, [t_1, v]) + \|\Delta y(\sigma)\| + \mathrm{Osc}(y, [\sigma, t_2])$$

$$\leq \left(1 + \frac{4}{a}\right)(D(\overline{\Pi}_m(\eta)) + 2\eta) + \mathbf{I}\eta + \Pi_{m-1}(\eta)$$

$$\leq \Pi_m(\mathrm{Osc}(x, [t_1, t_2]) + \delta).$$

On combining all of the cases above, we have

(31) $\qquad \mathrm{Osc}(w, [t_1, t_2]) \leq \Pi_m(\mathrm{Osc}(x, [t_1, t_2]) + \delta),$

(32) $\qquad \mathrm{Osc}(y, [t_1, t_2]) \leq \Pi_m(\mathrm{Osc}(x, [t_1, t_2]) + \delta).$

This completes the induction step. $\square$

REMARK. The proof of the above theorem was inspired by the proof of Lemma 4.3 of [4]. Because of the condition (i) in Theorem 4.1, the oscillation inequality given here is localized. Similar, but nonlocalized, oscillation inequalities were proved in [15] when $\overline{G} = \mathbb{R}_+^d$ and in [3] for a sequence of convex polyhedrons; in these cases, the direction of reflection was constant on each boundary face.

4.2. *C-tightness result.* Throughout this subsection and the next, we suppose that the following assumption holds in addition to (A1)–(A5).

ASSUMPTION 4.1. *There is a sequence of strictly positive constants* $\{\delta^n\}_{n=1}^{\infty}$ *such that for each positive integer $n$, there are processes $W^n, \widetilde{W}^n, X^n, \alpha^n$ having paths in $D([0,\infty), \mathbb{R}^d)$ and processes $Y^n, \widetilde{Y}^n, \beta^n$ having paths in $D([0,\infty), \mathbb{R}^\mathbf{I})$ defined on some probability space $(\Omega^n, \mathcal{F}^n, P^n)$ such that:*

(i) $P^n$-*a.s.*, $W^n = \widetilde{W}^n + \alpha^n$ *and* $\widetilde{W}^n(t) \in U_{\delta^n}(G)$ *for all* $t \geq 0$,

(ii) $P^n$-*a.s.*, $W^n(t) = X^n(t) + \sum_{i \in \mathcal{I}} \int_{(0,t]} \gamma^{i,n}(W^n(s-), W^n(s)) \, dY_i^n(s)$ *for all $t \geq 0$, where for each $i \in \mathcal{I}$, $\gamma^{i,n} : \mathbb{R}^d \times \mathbb{R}^d \to \mathbb{R}^d$ is Borel measurable and $\|\gamma^{i,n}(y, x)\| = 1$ for all $x, y \in \mathbb{R}^d$,*

(iii) $Y^n = \widetilde{Y}^n + \beta^n$, *where $\beta^n$ is locally of bounded variation and $P^n$-a.s., for each $i \in \mathcal{I}$,*

(a) $\widetilde{Y}_i^n(0) = 0$,
(b) $\widetilde{Y}_i^n$ *is nondecreasing and* $\Delta \widetilde{Y}_i^n(t) \leq \delta^n$ *for all* $t > 0$,
(c) $\widetilde{Y}_i^n(t) = \int_{(0,t]} 1_{\{\widetilde{W}^n(s) \in U_{\delta^n}(\partial G_i \cap \partial G)\}} \, d\widetilde{Y}_i^n(s),$

(iv) $\delta^n \to 0$ *as $n \to \infty$, and for each $\varepsilon > 0$, there is $\eta_\varepsilon > 0$ and $n_\varepsilon > 0$ such that for each $i \in \mathcal{I}$, $\|\gamma^{i,n}(y, x) - \gamma^i(x)\| < \varepsilon$ whenever $\|x - y\| < \eta_\varepsilon$ and $n \geq n_\varepsilon$,*

(v) $\alpha^n \to \mathbf{0}$ and $\mathcal{V}(\beta^n) \to \mathbf{0}$ in probability, as $n \to \infty$,
(vi) $\{X^n\}$ is $C$-tight.

REMARK. A simple case in which (iv) above holds is where $\gamma^{i,n}(y, x) \equiv \gamma^i(y)$. In (v), $\mathcal{V}(\beta^n)$ is the total variation process for $\beta^n$ (cf. Section 1.1).

The following theorem will play an important role in the proof of the invariance principle. It will be used to show that a sequence of processes satisfying suitably perturbed versions of the defining conditions for an SRBM [cf. (i)–(vi) above] is $C$-tight.

THEOREM 4.2 ($C$-tightness). *Suppose that Assumption 4.1 holds. Define $\mathcal{Z}^n = (W^n, X^n, Y^n)$ for each $n$. Then the sequence of processes $\{\mathcal{Z}^n\}_{n=1}^\infty$ is $C$-tight.*

REMARK. Note that $C$-tightness of $\{W^n\}$, $\{X^n\}$ and $\{Y^n\}$ implies $C$-tightness of $\{\mathcal{Z}^n\}$ (see Chapter VI, Corollary 3.33 of [12] for details).

PROOF OF THEOREM 4.2. References here to (i)–(vi) are to the conditions in Assumption 4.1.

Simple algebraic manipulations yield $P^n$-a.s.,

$$\widetilde{W}^n(t) = \widetilde{X}^n(t) + \sum_{i \in \mathcal{I}} \int_{(0,t]} \gamma^{i,n}(W^n(s-), W^n(s)) \, d\widetilde{Y}_i^n(s) \tag{33}$$

$$= \widetilde{X}^n(t) + \widetilde{V}^n(t) + \sum_{i \in \mathcal{I}} \int_{(0,t]} \gamma^i(\widetilde{W}^n(s)) \, d\widetilde{Y}_i^n(s), \tag{34}$$

where

$$\widetilde{X}^n(t) = X^n(t) + \left(-\alpha^n(t) + \sum_{i \in \mathcal{I}} \int_{(0,t]} \gamma^{i,n}(W^n(s-), W^n(s)) \, d\beta_i^n(s)\right) \tag{35}$$

and

$$\widetilde{V}^n(t) = \sum_{i \in \mathcal{I}} \int_{(0,t]} (\gamma^{i,n}(W^n(s-), W^n(s)) - \gamma^i(W^n(s))) \, d\widetilde{Y}_i^n(s) \tag{36}$$
$$+ \sum_{i \in \mathcal{I}} \int_{(0,t]} (\gamma^i(W^n(s)) - \gamma^i(\widetilde{W}^n(s))) \, d\widetilde{Y}_i^n(s).$$

The hypotheses on $\alpha^n$, the total variation process $\mathcal{V}(\beta^n)$ of $\beta^n$, and the fact that $\|\gamma^{i,n}(y, x)\| = 1$ for all $x, y \in \mathbb{R}^d$ and each $i \in \mathcal{I}$, imply that the process

$$-\alpha^n(\cdot) + \sum_{i \in \mathcal{I}} \int_{(0,\cdot]} \gamma^{i,n}(W^n(s-), W^n(s)) \, d\beta_i^n(s)$$

converges to **0** in probability as $n \to \infty$. Combining this with the fact that $\{X^n\}_{n=1}^\infty$ is $C$-tight, we obtain that $\{\widetilde{X}^n\}_{n=1}^\infty$ is $C$-tight.

Recall the positive nondecreasing function $\Pi(\cdot)$ from Theorem 4.1, and the constants $a$, $L$ and functions $R(\cdot)$ and $D(\cdot)$ from assumptions (A1)–(A5) in Section 3. Recall also that $\rho_0 = \frac{a}{4L}$.

Fix $\rho, \varepsilon, \eta, T$ such that $0 < \rho < \min\{\frac{\rho_0}{4}, \frac{R(a/4)}{4}\}$, $\varepsilon > 0$, $\eta > 0$ and $T > 0$. By assumption (A3), there is a constant $r_1 > 0$ such that

$$(37) \qquad D(r) < \min\left\{\frac{\rho}{2}, \varepsilon\right\} \qquad \text{for all } r \in (0, r_1].$$

Since $\Pi(u) \to 0$ as $u \to 0$, there are constants $0 < r_3 < r_2 < \min\{r_1, \frac{\rho}{4}, \frac{\varepsilon}{2\mathbf{I}}\}$ such that

$$(38) \qquad \Pi(r) < \frac{r_2}{2} \qquad \text{for all } r \in (0, r_3].$$

By (iv), there are $0 < \tilde{\varepsilon} < \min\{\frac{\rho}{8}, \frac{\varepsilon}{8}, \frac{r_3}{3}\}$ and $n_0 > 0$ such that for all $n \geq n_0$,

$$(39) \qquad \sup_{x \in \mathbb{R}^d} \sup_{\|y-x\| < 2\tilde{\varepsilon}} \max_{i \in \mathcal{I}} \|\gamma^{i,n}(y, x) - \gamma^i(x)\| < \frac{r_3}{6\mathbf{I}r_2}.$$

By (iv)–(vi), and Proposition 1.1, there exist an integer $n_1 > n_0$, a constant $\widetilde{M}_{\eta,T} > 0$ and $\widetilde{\lambda} \in (0, T)$, such that for all $n \geq n_1$,

$$(40) \qquad P^n\left\{\sup_{0 \leq s \leq T} \|\widetilde{X}^n(s)\| \leq \widetilde{M}_{\eta,T}\right\} \geq 1 - \eta/2,$$

$$(41) \qquad P^n\{w_T(\widetilde{X}^n, \widetilde{\lambda}) \geq \tilde{\varepsilon}\} \leq \eta/4,$$

$$(42) \qquad P^n\left\{\sup_{0 \leq s \leq T} \|\alpha^n(s)\| < \frac{r_3}{6\mathbf{I}Lr_2} \wedge \frac{\tilde{\varepsilon}}{8}\right\} \geq 1 - \eta/4,$$

$$(43) \qquad \delta^n < \min\left\{\frac{r_3}{3}, \frac{r_2}{2\mathbf{I}}, \frac{\rho}{8(1+\mathbf{I})}, \frac{\tilde{\varepsilon}}{2\mathbf{I}}\right\}.$$

To prove $C$-tightness of $\{\widetilde{W}^n\}$ and $\{\widetilde{Y}^n\}$ (and hence of $\{W^n\}$, $\{Y^n\}$), by Proposition 1.1, it suffices to show that there exists a constant $N_{\eta,T} > 0$ such that for all $n \geq n_1$,

$$(44) \qquad P^n\{w_T(\widetilde{W}^n, \widetilde{\lambda}) \geq \varepsilon\} \leq \eta,$$

$$(45) \qquad P^n\{w_T(\widetilde{Y}^n, \widetilde{\lambda}) \geq \varepsilon\} \leq \eta,$$

$$(46) \qquad P^n\left\{\sup_{0 \leq s \leq T} \|\widetilde{W}^n(s)\| \leq N_{\eta,T}\right\} \geq 1 - \eta,$$

$$(47) \qquad P^n\left\{\sup_{0 \leq s \leq T} \|\widetilde{Y}^n(s)\| \leq N_{\eta,T}\right\} \geq 1 - \eta.$$

For each $n \geq 1$, let $F^n$ be a set in $\mathcal{F}^n$ such that $P^n(F^n) = 1$ and on $F^n$, properties (iii)(a)–(c) hold, (33)–(36) hold, and $\widetilde{W}^n(t) \in U_{\delta^n}(G)$ for all $t \geq 0$. Fix a $t$ such that $0 \leq t < t + \widetilde{\lambda} \leq T$. Let

$$\tau^n = \inf\{s \geq t : \widetilde{W}^n(s) \in U_{\delta^n}(\partial G_i \cap \partial G) \text{ for some } i \in \mathcal{I}\}. \tag{48}$$

For each $n \geq n_1$, let

$$H^n = \left\{ w_T(\widetilde{X}^n, \widetilde{\lambda}) < \widetilde{\varepsilon}, \sup_{0 \leq s \leq T} \|\alpha^n(s)\| < \frac{r_3}{6\mathbf{I}Lr_2} \wedge \frac{\widetilde{\varepsilon}}{8}, \right.$$
$$\left. \sup_{0 \leq s \leq T} \|\widetilde{X}^n(s)\| \leq \widetilde{M}_{\eta,T} \right\} \cap F^n. \tag{49}$$

Then by (40)–(42) and the definition of $F^n$,

$$P\{H^n\} \geq 1 - \eta. \tag{50}$$

Fix $\omega^n \in H^n$. By the definition of $w_T(x, \lambda)$ in (3), we have that,

$$\sup_{r,s \in [t,t+\widetilde{\lambda}]} \|\widetilde{X}^n(s, \omega^n) - \widetilde{X}^n(r, \omega^n)\| < \widetilde{\varepsilon}. \tag{51}$$

Now there are two cases to consider for $n \geq n_1$ and $u$, $v$ fixed such that $t \leq u < v \leq t + \widetilde{\lambda}$.

*Case 1.* $\omega^n \in \{\tau^n > v\}$. In this case, by (iii)(c), $\widetilde{Y}^n(\cdot, \omega^n)$ does not increase on the interval $(u, v]$, that is, $\widetilde{Y}_i^n(v, \omega^n) - \widetilde{Y}_i^n(u, \omega^n) = 0$ for all $i \in \mathcal{I}$. Then by (34) and (36),

$$\widetilde{W}^n(v, \omega^n) - \widetilde{W}^n(u, \omega^n) = \widetilde{X}^n(v, \omega^n) - \widetilde{X}^n(u, \omega^n). \tag{52}$$

Hence, by (51),

$$\|\widetilde{W}^n(v, \omega^n) - \widetilde{W}^n(u, \omega^n)\| \leq \sup_{r,s \in [t,t+\widetilde{\lambda}]} \|\widetilde{X}^n(s, \omega^n) - \widetilde{X}^n(r, \omega^n)\| < \widetilde{\varepsilon} < \varepsilon/8,$$

and we also have

$$\|\widetilde{Y}^n(v, \omega^n) - \widetilde{Y}^n(u, \omega^n)\| = 0 < \varepsilon/2.$$

*Case 2.* $\omega^n \in \{\tau^n \leq v\}$. Then there is an $i \in \mathcal{I}$ such that $\widetilde{W}^n(\tau^n, \omega^n) \in U_{\delta^n}(\partial G_i \cap \partial G)$, since the set $U_{\delta^n}(\partial G_i \cap \partial G)$ is closed and $\widetilde{W}^n(\cdot, \omega^n)$ is right continuous. It follows that there is some $x_0 \in \partial G$ (which depends on $\omega^n$) such that $\widetilde{W}^n(\tau^n, \omega^n)$ is in the closed ball $B_{\delta^n}(x_0) \subset B_\rho(x_0)$. To apply the oscillation inequality in Theorem 4.1, we first prove the following:

(53) $\widetilde{W}^n(r,\omega^n) \in B_\rho(x_0)$ for all $r$ satisfying $\tau^n \leq r \leq v$.

For the proof of (53), let

(54) $$\xi^n = \inf\{r \geq \tau^n : \widetilde{W}^n(r,\omega^n) \notin B_\rho(x_0)\} \wedge v.$$

By the definition of $\xi^n$, $\widetilde{W}^n(r,\omega^n) \in B_\rho(x_0)$ for each $r \in [\tau^n, \xi^n)$. In order to apply the oscillation inequality in Theorem 4.1 on the time interval $[\tau^n, \xi^n)$, we show that

(55) $$D(\Pi(\operatorname{Osc}(\widetilde{X}^n(\cdot,\omega^n) + \widetilde{V}^n(\cdot,\omega^n), [\tau^n, \xi^n)) + \delta^n)) < \frac{\rho}{2}.$$

For each $r \in (0,T]$, by (i)–(iii) and (33), (49), (43), we have that

$$\|W^n(r-,\omega^n) - W^n(r,\omega^n)\|$$
$$\leq \|\widetilde{W}^n(r-,\omega^n) - \widetilde{W}^n(r,\omega^n)\| + \|\alpha^n(r-,\omega^n) - \alpha^n(r,\omega^n)\|$$
$$\leq \|\Delta \widetilde{X}^n(r,\omega^n)\| + 2\sup_{0 \leq s \leq T}\|\alpha^n(s)\| + \mathbf{I}\delta^n$$
$$\leq \tilde{\varepsilon} + \frac{\tilde{\varepsilon}}{4} + \frac{\tilde{\varepsilon}}{2} < 2\tilde{\varepsilon}.$$

Hence by (39), for each $r \in (0,T]$,

(56) $$\|\gamma^{i,n}(W^n(r-,\omega^n), W^n(r,\omega^n)) - \gamma^i(W^n(r,\omega^n))\| \leq \frac{r_3}{6\mathbf{I}r_2}.$$

By (36), (56), Assumption (A4), (i) and (49), we have that for any $s_1$, $s_2$ such that $u \leq s_1 < s_2 \leq v$,

(57)
$$\|\widetilde{V}^n(s_2,\omega^n) - \widetilde{V}^n(s_1,\omega^n)\|$$
$$\leq \sum_{i \in \mathcal{I}} \int_{(s_1,s_2]} \|\gamma^{i,n}(W^n(r-,\omega^n), W^n(r,\omega^n))$$
$$\qquad\qquad - \gamma^i(W^n(r,\omega^n))\| \, d\widetilde{Y}_i^n(r,\omega^n)$$
$$+ \sum_{i \in \mathcal{I}} \int_{(s_1,s_2]} \|\gamma^i(W^n(r,\omega^n)) - \gamma^i(\widetilde{W}^n(r,\omega^n))\| \, d\widetilde{Y}_i^n(r,\omega^n)$$
$$\leq \sum_{i \in \mathcal{I}} \frac{r_3}{6\mathbf{I}r_2}(\widetilde{Y}_i^n(s_2,\omega^n) - \widetilde{Y}_i^n(s_1,\omega^n))$$
$$+ \sum_{i \in \mathcal{I}} \int_{(s_1,s_2]} L\|W^n(r,\omega^n) - \widetilde{W}^n(r,\omega^n)\| \, d\widetilde{Y}_i^n(r,\omega^n)$$
$$\leq \sum_{i \in \mathcal{I}} \frac{r_3}{6\mathbf{I}r_2}(\widetilde{Y}_i^n(s_2,\omega^n) - \widetilde{Y}_i^n(s_1,\omega^n))$$
$$+ \sum_{i \in \mathcal{I}} L\frac{r_3}{6\mathbf{I}Lr_2}(\widetilde{Y}_i^n(s_2,\omega^n) - \widetilde{Y}_i^n(s_1,\omega^n))$$

$$\leq \frac{r_3}{3r_2}\|\widetilde{Y}^n(s_2,\omega^n) - \widetilde{Y}^n(s_1,\omega^n)\|.$$

Let

(58) $$\sigma^n = \inf\{s \geq \tau^n : \operatorname{Osc}(\widetilde{Y}^n(\cdot,\omega^n),[\tau^n,s)) > r_2\}.$$

Note that $\operatorname{Osc}(\widetilde{Y}^n(\cdot,\omega^n),[\tau^n,s))$ as a function of $s$ defined on $(\tau^n,\infty)$ is left continuous with finite right limits and is nondecreasing. By the right continuity of $\widetilde{Y}^n$, we know that

$$\operatorname{Osc}(\widetilde{Y}^n(\cdot,\omega^n),[\tau^n,s)) \to 0 \qquad \text{as } s \downarrow \tau^n.$$

Thus, $\sigma^n > \tau^n$, $\operatorname{Osc}(\widetilde{Y}^n(\cdot,\omega^n),[\tau^n,\sigma^n)) \leq r_2$ and on $\{\sigma^n < \infty\}$, $\operatorname{Osc}(\widetilde{Y}^n(\cdot,\omega^n),[\tau^n,\sigma^n]) \geq r_2$. By (57), (51), (43), the choice of $\varepsilon$, and since $t \leq \tau^n \leq \xi^n \leq v \leq t + \widetilde{\lambda}$, we have

(59)
$$\begin{aligned}
\operatorname{Osc}&(\widetilde{X}^n(\cdot,\omega^n) + \widetilde{V}^n(\cdot,\omega^n),[\tau^n,\xi^n \wedge \sigma^n)) + \delta^n \\
&\leq \operatorname{Osc}(\widetilde{X}^n(\cdot,\omega^n),[\tau^n,\xi^n \wedge \sigma^n)) \\
&\quad + \operatorname{Osc}(\widetilde{V}^n(\cdot,\omega^n),[\tau^n,\xi^n \wedge \sigma^n)) + \delta^n \\
&\leq \operatorname{Osc}(\widetilde{X}^n(\cdot,\omega^n),[\tau^n,\xi^n \wedge \sigma^n)) \\
&\quad + \frac{r_3}{3r_2}\operatorname{Osc}(\widetilde{Y}^n(\cdot,\omega^n),[\tau^n,\xi^n \wedge \sigma^n)) + \delta^n \\
&\leq \widetilde{\varepsilon} + \frac{r_3}{3r_2}r_2 + \delta^n < r_3.
\end{aligned}$$

Then by (38) and the monotonicity of $D(\cdot)$, we have

(60)
$$\begin{aligned}
D(\Pi(\operatorname{Osc}(\widetilde{X}^n(\cdot,\omega^n) &+ \widetilde{V}^n(\cdot,\omega^n),[\tau^n,\xi^n \wedge \sigma^n)) + \delta^n)) \\
&\leq D\left(\frac{r_2}{2}\right) \leq D(r_2) \leq D(r_1) < \frac{\rho}{2}.
\end{aligned}$$

We claim that

(61) $$\sigma^n \geq \xi^n.$$

To prove (61), we proceed by contradiction and suppose that $\sigma^n < \xi^n$. Then by (60), with $x = \widetilde{X}^n(\cdot,\omega^n) + \widetilde{V}^n(\cdot,\omega^n)$ and $\delta = \delta^n$, condition (iv) of Theorem 4.1 holds with $[s,t] = [\tau^n,\sigma^n - 1/m]$ for all $m$ sufficiently large. By applying Theorem 4.1 and letting $m \to \infty$, we obtain using (34), (38) and (59) that,

(62)
$$\begin{aligned}
\operatorname{Osc}&(\widetilde{Y}^n(\cdot,\omega^n),[\tau^n,\sigma^n)) \\
&\leq \Pi(\operatorname{Osc}(\widetilde{X}^n(\cdot,\omega^n) + \widetilde{V}^n(\cdot,\omega^n),[\tau^n,\xi^n \wedge \sigma^n)) + \delta^n) \\
&\leq \Pi(r_3) < \frac{r_2}{2}.
\end{aligned}$$

By (62), (iii)(b) and (43), we obtain that

$$\mathrm{Osc}(\widetilde{Y}^n(\cdot,\omega^n),[\tau^n,\sigma^n]) \leq \frac{r_2}{2} + \mathbf{I}\delta^n < r_2.$$

This contradicts the fact that $\mathrm{Osc}(\widetilde{Y}^n(\cdot,\omega^n),[\tau^n,\sigma^n]) \geq r_2$ on $\{\sigma^n < \infty\}$, and so (61) holds and (55) follows by (60).

By applying Theorem 4.1 on $[\tau^n, \xi^n - 1/m]$ and then letting $m \to \infty$, we obtain using (61), (59) and (38), that

$$\begin{aligned}\mathrm{Osc}(\widetilde{W}^n(\cdot,\omega^n),[\tau^n,\xi^n)) \\ &\leq \Pi(\mathrm{Osc}(\widetilde{X}^n(\cdot,\omega^n)+\widetilde{V}^n(\cdot,\omega^n),[\tau^n,\xi^n\wedge\sigma^n))+\delta^n) \\ &< \frac{r_2}{2},\end{aligned}$$

and similarly,

(63) $$\mathrm{Osc}(\widetilde{Y}^n(\cdot,\omega^n),[\tau^n,\xi^n)) < \frac{r_2}{2}.$$

Then we have

$$\begin{aligned}\|\widetilde{W}^n(\xi^n-,\omega^n) - x_0\| \\ &\leq \|\widetilde{W}^n(\xi^n-,\omega^n) - \widetilde{W}^n(\tau^n,\omega^n)\| + \|\widetilde{W}^n(\tau^n,\omega^n) - x_0\| \\ &\leq \frac{r_2}{2} + \delta^n.\end{aligned}$$

Using hypotheses (ii), (iii)(b), and (33), (51), we obtain

$$\begin{aligned}\|\widetilde{W}^n(\xi^n,\omega^n) - \widetilde{W}^n(\xi^n-,\omega^n)\| \\ &\leq \|\widetilde{X}^n(\xi^n,\omega^n) - \widetilde{X}^n(\xi^n-,\omega^n)\| \\ &\quad + \sum_{i\in\mathcal{I}}\|\gamma^{i,n}(W^n(\xi^n-,\omega^n),W^n(\xi^n,\omega^n))\| \\ &\qquad\qquad \times (\widetilde{Y}^n_i(\xi^n,\omega^n) - \widetilde{Y}^n_i(\xi^n-,\omega^n)) \\ &\leq \tilde{\varepsilon} + \mathbf{I}\delta^n.\end{aligned}$$

Hence

$$\begin{aligned}\|\widetilde{W}^n(\xi^n,\omega^n) - x_0\| &\leq \|\widetilde{W}^n(\xi^n-,\omega^n) - x_0\| \\ &\quad + \|\widetilde{W}^n(\xi^n,\omega^n) - \widetilde{W}^n(\xi^n-,\omega^n)\| \\ &\leq \frac{r_2}{2} + \delta^n + \tilde{\varepsilon} + \mathbf{I}\delta^n \\ &\leq \tilde{\varepsilon} + (\mathbf{I}+1)\delta^n + \frac{r_2}{2} \\ &< \rho/8 + \rho/8 + \rho/8 < \rho/2.\end{aligned}$$

It follows from this that $\xi^n = v$ and (53) holds, as desired.

Then, by (33), (51), (iii)(b), (iii)(c), (63) and (43), we have

$$\|\widetilde{W}^n(v,\omega^n) - \widetilde{W}^n(u,\omega^n)\|$$
$$\leq \sup_{r,s \in [u,v]} \|\widetilde{X}^n(s,\omega^n) - \widetilde{X}^n(r,\omega^n)\| + \sum_{i \in \mathcal{I}} (\widetilde{Y}_i^n(v,\omega^n) - \widetilde{Y}_i^n(u,\omega^n))$$
$$\leq \tilde{\varepsilon} + \sum_{i \in \mathcal{I}} (\widetilde{Y}_i^n(v,\omega^n) - \widetilde{Y}_i^n(u \vee \tau^n, \omega^n))$$
(64)
$$+ \sum_{i \in \mathcal{I}} (\widetilde{Y}_i^n(u \vee \tau^n, \omega^n) - \widetilde{Y}_i^n(u,\omega^n))$$
$$\leq \tilde{\varepsilon} + \mathbf{I} \operatorname{Osc}(\widetilde{Y}^n(\cdot,\omega^n), [u \vee \tau^n, v)) + \sum_{i \in \mathcal{I}} \Delta \widetilde{Y}_i^n(v,\omega^n) + \mathbf{I}\delta^n$$
$$\leq \tilde{\varepsilon} + \mathbf{I}\frac{r_2}{2} + \mathbf{I}\delta^n + \mathbf{I}\delta^n < \frac{\varepsilon}{8} + \frac{\varepsilon}{4} + \frac{\varepsilon}{16} + \frac{\varepsilon}{16} = \frac{\varepsilon}{2}$$

and

(65)
$$\|\widetilde{Y}^n(v,\omega^n) - \widetilde{Y}^n(u,\omega^n)\| \leq \sum_{i \in \mathcal{I}} (\widetilde{Y}_i^n(v,\omega^n) - \widetilde{Y}_i^n(u,\omega^n))$$
$$\leq \sum_{i \in \mathcal{I}} (\widetilde{Y}_i^n(v,\omega^n) - \widetilde{Y}_i^n(u \vee \tau^n, \omega^n))$$
$$+ \sum_{i \in \mathcal{I}} (\widetilde{Y}_i^n(u \vee \tau^n, \omega^n) - \widetilde{Y}_i^n(u,\omega^n))$$
$$< \frac{\varepsilon}{2}.$$

Here we have used the fact that $\tilde{Y}_i$ does not increase on $(u, \tau^n \vee u)$ and can jump at most by $\delta^n$ at $\tau^n$, by the definition of $\tau^n$ and (iii)(c).

On combining the results from Case 1 and Case 2, we obtain that for each $n \geq n_1$,

(66) $\sup \left\{ \sup_{u,v \in [t, t+\widetilde{\lambda}]} \|\widetilde{W}^n(v,\omega^n) - \widetilde{W}^n(u,\omega^n)\| : 0 \leq t \leq t + \widetilde{\lambda} \leq T \right\} \leq \frac{\varepsilon}{2} < \varepsilon$

and

(67) $\sup \left\{ \sup_{u,v \in [t, t+\widetilde{\lambda}]} \|\widetilde{Y}^n(v,\omega^n) - \widetilde{Y}^n(u,\omega^n)\| : 0 \leq t \leq t + \widetilde{\lambda} \leq T \right\} \leq \frac{\varepsilon}{2} < \varepsilon.$

Hence since $\omega^n \in H^n$ was arbitrary, by (50), we have that (44) and (45) hold for all $n \geq n_1$.

Next we show that there is a constant $N_{\eta,T} > 0$ such that (46) and (47) hold for all $n \geq n_1$. By (66)–(67) above, we have that for each $n \geq n_1$, $\omega^n \in$ $H^n$, $t$ such that $0 \leq t < t + \widetilde{\lambda} \leq T$ and $t \leq u < v \leq t + \widetilde{\lambda}$,

(68) $$\|\widetilde{W}^n(v,\omega^n) - \widetilde{W}^n(u,\omega^n)\| < \varepsilon$$

and

(69) $$\|\widetilde{Y}^n(v,\omega^n) - \widetilde{Y}^n(u,\omega^n)\| < \varepsilon.$$

Then, for each $0 \leq s \leq T$, by (68), (69), (49) and (33), we have

$$\|\widetilde{W}^n(s,\omega^n)\| \leq \|\widetilde{W}^n(s,\omega^n) - \widetilde{W}^n(0,\omega^n)\| + \|\widetilde{W}^n(0,\omega^n)\|$$
$$\leq \sum_{i=1}^{[T/\widetilde{\lambda}]+1} \|\widetilde{W}^n(i\widetilde{\lambda} \wedge s, \omega^n) - \widetilde{W}^n((i-1)\widetilde{\lambda} \wedge s, \omega^n)\| + \|\widetilde{X}^n(0,\omega^n)\|$$
$$\leq ([T/\widetilde{\lambda}] + 1)\varepsilon + \widetilde{M}_{\eta,T}$$

and

$$\|\widetilde{Y}^n(s,\omega^n)\| \leq \|\widetilde{Y}^n(s,\omega^n) - \widetilde{Y}^n(0,\omega^n)\|$$
$$\leq \sum_{i=1}^{[T/\widetilde{\lambda}]+1} \|\widetilde{Y}^n(i\widetilde{\lambda} \wedge s, \omega^n) - \widetilde{Y}^n((i-1)\widetilde{\lambda} \wedge s, \omega^n)\|$$
$$\leq ([T/\widetilde{\lambda}] + 1)\varepsilon.$$

Here $[T/\widetilde{\lambda}]$ is the greatest integer less than or equal to $T/\widetilde{\lambda}$. Let $N_{\eta,T} = ([T/\widetilde{\lambda}] + 1)\varepsilon + \widetilde{M}_{\eta,T}$. Then we obtain that for $n \geq n_1$ and $\omega^n \in H^n$,

(70) $$\sup_{0 \leq s \leq T} \|\widetilde{W}^n(s,\omega^n)\| \leq N_{\eta,T}$$

and

(71) $$\sup_{0 \leq s \leq T} \|\widetilde{Y}^n(s,\omega^n)\| \leq N_{\eta,T}.$$

Then by (50), we have that (46) and (47) hold for all $n \geq n_1$.

Finally by applying Proposition 1.1, we have the $C$-tightness of $\{\widetilde{W}^n\}$ and $\{\widetilde{Y}^n\}$. It then follows that $\{(\widetilde{W}^n, X^n, \widetilde{Y}^n)\}_{n=1}^{\infty}$ is $C$-tight. Since $\mathcal{Z}^n = (\widetilde{W}^n, X^n, \widetilde{Y}^n) + (\alpha^n, \mathbf{0}, \beta^n)$ where $\alpha^n, \mathcal{V}(\beta^n) \to \mathbf{0}$ in probability as $n \to \infty$, then $\{\mathcal{Z}^n\}_{n=1}^{\infty}$ is also $C$-tight. □

4.3. *Invariance principle for SRBMs.* The main theorem of the paper is the following.

THEOREM 4.3 (Invariance principle for SRBMs). *Suppose that Assumption* 4.1 *holds. Define* $\mathcal{Z}^n = (W^n, X^n, Y^n)$ *for each* $n$. *Then the sequence of* processes $\{\mathcal{Z}^n\}_{n=1}^\infty$ is C-tight and any (weak) limit point of this sequence is of the form $\mathcal{Z} = (W, X, Y)$ where continuous d-dimensional processes $W, X$ and a continuous **I**-dimensional process $Y$ are defined on some probability space $(\Omega, \mathcal{F}, P)$ such that conditions (i), (ii) and (iv) of Definition 2.1 hold with $\mathcal{F}_t = \sigma\{\mathcal{Z}(s) : 0 \leq s \leq t\}$, $t \geq 0$.

If, in addition, the following conditions (vi)′ and (vii) hold, then any weak limit point of the sequence $\{\mathcal{Z}^n\}_{n=1}^\infty$ is an extended SRBM associated with the data $(G, \mu, \Gamma, \{\gamma^i, i \in \mathcal{I}\}, \nu)$. If furthermore the following condition (viii) holds, then $W^n \Rightarrow W$ as $n \to \infty$ where $W$ is an SRBM associated with $(G, \mu, \Gamma, \{\gamma^i, i \in \mathcal{I}\}, \nu)$.

(vi)′ $\{X^n\}$ converges in distribution to a d-dimensional Brownian motion with drift $\mu$, covariance matrix $\Gamma$ and initial distribution $\nu$.

(vii) For each (weak) limit point $\mathcal{Z} = (W, X, Y)$ of $\{\mathcal{Z}^n\}_{n=1}^\infty$, $\{X(t) - X(0) - \mu t, \mathcal{F}_t, t \geq 0\}$ is a martingale.

(viii) If a process $W$ satisfies the properties in Definition 2.1, the law of $W$ is unique, that is, the law of an SRBM associated with the data $(G, \mu, \Gamma, \{\gamma^i, i \in \mathcal{I}\}, \nu)$ is unique.

REMARK. We note that (vi)′ implies that (vi) of Assumption 4.1 holds.

PROOF OF THEOREM 4.3. By Theorem 4.2, we have that the sequence $\{\mathcal{Z}^n\}_{n=1}^\infty$ is C-tight. Let $\mathcal{Z} = (W, X, Y)$ be a (weak) limit point of $\{\mathcal{Z}^n\}_{n=1}^\infty$, that is, there is a subsequence $\{n_k\}$ of $\{n\}$ such that $\mathcal{Z}^{n_k} \Rightarrow \mathcal{Z}$ as $k \to \infty$. It also follows that $\widetilde{\mathcal{Z}}^{n_k} \equiv (\widetilde{W}^{n_k}, X^{n_k}, \widetilde{Y}^{n_k}) \Rightarrow \mathcal{Z}$ as $k \to \infty$. By the C-tightness of $\{\mathcal{Z}^n\}$, we obtain that $\mathcal{Z}$ has continuous paths a.s. For the purpose of verifying that $\mathcal{Z}$ satisfies the listed properties in Definition 2.1, one may invoke the Skorokhod representation theorem to assume, without loss of generality, that $\mathcal{Z}^{n_k}$ and $\widetilde{\mathcal{Z}}^{n_k}$ converge u.o.c. to $\mathcal{Z}$ a.s. as $k \to \infty$ and $\mathcal{V}(\beta^{n_k})$ converges u.o.c. to $\mathbf{0}$ a.s. as $k \to \infty$. With this simplification, it is easily verified that the properties of $\{\mathcal{Z}^{n_k}\}$ and $\{\widetilde{\mathcal{Z}}^{n_k}\}$ imply that $\mathcal{Z}$ has properties (ii) and (iv)(a)–(b) of Definition 2.1. For the verification of property (i) of Definition 2.1, note that for each $k$, a.s. for each $t \geq 0$,

$$
\begin{aligned}
W^{n_k}(t) = {} & X^{n_k}(t) + \sum_{i \in \mathcal{I}} \int_{(0,t]} \gamma^{i,n_k}(W^{n_k}(s-), W^{n_k}(s)) \, d\beta_i^{n_k}(s) \\
& + \sum_{i \in \mathcal{I}} \int_{(0,t]} (\gamma^{i,n_k}(W^{n_k}(s-), W^{n_k}(s)) - \gamma^i(W^{n_k}(s))) \, d\widetilde{Y}_i^{n_k}(s) \\
& + \sum_{i \in \mathcal{I}} \int_{(0,t]} \gamma^i(W^{n_k}(s)) \, d\widetilde{Y}_i^{n_k}(s).
\end{aligned}
$$

The sum of the first two terms on the right-hand side of the above equality converges a.s. to $X(t)$ as $k \to \infty$. The third term on the right-hand side converges a.s. to $\mathbf{0}$ as $k \to \infty$, by property (iv) and the fact that a.s.,

$$\sup_{s \in (0,t]} \|W^{n_k}(s) - W^{n_k}(s-)\|$$
$$\leq \sup_{s \in (0,t]} \|\Delta X^{n_k}(s)\| + \mathbf{I} \sup_{s \in (0,t]} \|\Delta Y^{n_k}(s)\| \to 0 \qquad \text{as } k \to \infty.$$

It remains to show that for each $i \in \mathcal{I}$ and $t \geq 0$, a.s.,

$$\int_{(0,t]} \gamma^i(W^{n_k}(s)) \, d\widetilde{Y}_i^{n_k}(s) \to \int_{(0,t]} \gamma^i(W(s)) \, dY_i(s) \qquad \text{as } k \to \infty.$$

This follows directly from Lemma A.4.

For the verification of property (iv)(c) of Definition 2.1, it suffices to show that for each $i \in \mathcal{I}$, $m = 1, 2, \ldots$, a.s. for each $t \geq 0$,

(72) $$Y_i(t) = \int_{(0,t]} f_m(W(s)) \, dY_i(s),$$

where $\{f_m\}_{m=1}^\infty$ is a sequence of real valued continuous functions defined on $\mathbb{R}^d$ such that for each $m$, the range of $f_m$ is $[0,1]$, $f_m(x) = 1$ for $x \in U_{1/m}(\partial G_i \cap \partial G)$ and $f_m(x) = 0$ for $x \notin U_{2/m}(\partial G_i \cap \partial G)$. The existence of such a sequence of continuous functions $\{f_m\}_{m=1}^\infty$ can be shown using Urysohn's lemma (cf. [8], page 122). Then (72) is a consequence of Lemma A.4, property (iii) of $\widetilde{Y}_i^{n_k}$ and the fact that $\delta^{n_k} \to 0$ as $k \to \infty$. Indeed, a.s., for each $t \geq 0$,

$$Y_i(t) = \lim_{k \to \infty} \widetilde{Y}_i^{n_k}(t) = \lim_{k \to \infty} \int_{(0,t]} 1_{\{\widetilde{W}^{n_k}(s) \in U_{\delta^{n_k}}(\partial G_i \cap \partial G)\}} \, d\widetilde{Y}_i^{n_k}(s)$$
$$= \lim_{k \to \infty} \int_{(0,t]} f_m(\widetilde{W}^{n_k}(s)) \, d\widetilde{Y}_i^{n_k}(s)$$
$$= \int_{(0,t]} f_m(W(s)) \, dY_i(s).$$

Thus, $\mathcal{Z}$ satisfies properties (i), (ii) and (iv) of Definition 2.1 with $\mathcal{F}_t = \sigma\{\mathcal{Z}(s) : 0 \leq s \leq t\}$, $t \geq 0$.

Assuming properties (vi)$'$ and (vii) hold, $\mathcal{Z}$ satisfies (iii) of Definition 2.1. Then $\mathcal{Z}$ is an extended SRBM associated with the data $(G, \mu, \Gamma, \{\gamma^i, i \in \mathcal{I}\}, \nu)$. If in addition, property (viii) holds, then the law of $W$ is unique. Since each weak limit $W$ is an SRBM associated with the data $(G, \mu, \Gamma, \{\gamma^i, i \in \mathcal{I}\}, \nu)$ and the law of such an SRBM is unique, then by a standard argument, $W^n \Rightarrow W$ as $n \to \infty$ where $W$ is an SRBM associated with $(G, \mu, \Gamma, \{\gamma^i, i \in \mathcal{I}\}, \nu)$. □

Some sufficient conditions for (vii) to hold are given in Proposition 4.2 of [15] for a simpler setting where $\overline{G} = \mathbb{R}^d_+$. Two of those conditions generalize to our setting here and can be proved in the same manner as in [15]. For completeness, we state the ensuing result here.

PROPOSITION 4.1. *Suppose that Assumption 4.1 and* (vi)′ *of Theorem 4.3 hold. If, in addition, one of the following conditions* (I)–(II) *holds, then condition* (vii) *of Theorem 4.3 is satisfied, and any weak limit point of* $\{\mathcal{Z}^n\}_{n=1}^\infty$ *is an extended SRBM associated with* $(G, \mu, \Gamma, \{\gamma^i, i \in \mathcal{I}\}, \nu)$.

(I) *For any triple of d-dimensional $\{\mathcal{F}_t\}$-adapted processes $(W, X, Y)$ defined on some filtered probability space $(\Omega, \mathcal{F}, \{\mathcal{F}_t\}, P)$ and satisfying conditions* (i), (ii) *and* (iv) *of Definition 2.1 together with the condition that $X$, under $P$, is a d-dimensional Brownian motion with drift vector $\mu$, covariance matrix $\Gamma$ and initial distribution $\nu$, the pair $(W, Y)$ is adapted to the filtration generated by $X$ and the $P$-null sets.*

(II) $X^n = \check{X}^n + \varepsilon_1^n$, $Y^n = \check{Y}^n + \varepsilon_2^n$, $W^n = \check{W}^n + \varepsilon_3^n$, *where $\varepsilon_1^n, \varepsilon_2^n, \varepsilon_3^n$ are processes converging to $\mathbf{0}$ in probability as $n \to \infty$, and:*

(a) $\{\check{X}^n(t) - \check{X}^n(0)\}_{n=1}^\infty$ *is uniformly integrable for each $t \geq 0$,*

(b) *there is a sequence of constants $\{\mu^n\}_{n=1}^\infty$ in $\mathbb{R}^d$ such that* $\lim_{n\to\infty} \mu^n = \mu$,

(c) *for each $n$, $\{\check{X}^n(t) - \check{X}^n(0) - \mu^n t, t \geq 0\}$ is a $P^n$-martingale with respect to the filtration generated by $(\check{W}^n, \check{X}^n, \check{Y}^n)$.*

In the rest of this work, we focus on applications of the invariance principle and in particular on giving sufficient conditions for property (viii) of Theorem 4.3 to hold.

**5. Applications of the invariance principle.** In Section 5.1, we prove weak existence of SRBMs associated with data $(G, \mu, \Gamma, \{\gamma^i, i \in \mathcal{I}\}, \nu)$ satisfying (A1)–(A5) of Section 3. This is accomplished by constructing a sequence of approximations whose weak limit points are SRBMs. The invariance principle is used to prove the $C$-tightness of the approximations and that any weak limit point is an SRBM. In Sections 5.2 and 5.3, using known results on uniqueness in law for SRBMs, we illustrate the invariance principle for certain domains and directions of reflection.

5.1. *Weak existence of SRBMs.*

THEOREM 5.1. *Suppose that assumptions* (A1)–(A5) *of Section 3 hold. Then there exists an SRBM associated with the data* $(G, \mu, \Gamma, \{\gamma^i, i \in \mathcal{I}\}, \nu)$.

PROOF. We construct a sequence of approximations to an SRBM and use the invariance principle to establish weak convergence along a subsequence to an SRBM.

In the following we will use $R(\cdot)$ from assumption (A2), $L > 0$ from assumption (A4), $a > 0$ from assumption (A5), and $\rho_0 = \frac{a}{4L}$. Fix $\varepsilon > 0$ and $0 < \rho < \min\{\frac{\rho_0}{4}, \frac{R(a/4)}{4}\}$. By assumption (A3), there is a constant $r_1 > 0$ such that

$$D(r) < \min\left\{\frac{\rho}{2}, \varepsilon\right\} \qquad \text{for all } r \in (0, r_1].$$

Recall the properties of $\Pi(\cdot)$ from Theorem 4.1. Since $\Pi(u) \to 0$ as $u \to 0$, there are constants $0 < r_3 < r_2 < \min\{r_1, \frac{\rho}{4}, \frac{\varepsilon}{2\mathbf{I}}\}$ such that

$$\Pi(r) < \frac{r_2}{2} \qquad \text{for all } r \in (0, r_3].$$

Fix $\tilde{\varepsilon}$ and $\delta$ such that $0 < \tilde{\varepsilon} < \min\{\frac{\rho}{8}, \frac{\varepsilon}{8}, \frac{r_3}{3}, \frac{r_3}{24\mathbf{I}Lr_2}\}$ and $0 < 2\delta < \min\{\frac{r_3}{3}, \frac{r_2}{2\mathbf{I}}, \frac{\rho}{8(1+\mathbf{I})}, \frac{\tilde{\varepsilon}}{2\mathbf{I}}\}$.

We will construct a $d$-dimensional stochastic process $W^\delta$ and an $\mathbf{I}$-dimensional "pushing" process $Y^\delta$, such that $W^\delta$ approximately satisfies the conditions defining an SRBM for the data $(G, \mu, \Gamma, \{\gamma^i, i \in \mathcal{I}\}, \nu)$ (cf. Assumption 4.1). The idea for this construction is to use a Brownian motion $X$ with drift vector $\mu$, covariance matrix $\Gamma$ and initial distribution $\nu$. Away from $\partial G$, the increments of $W^\delta$ are determined by those of $X$. For any time $t \geq 0$ such that $W^\delta(t-) \in \partial G$, we add an instantaneous jump to $W^\delta(t-)$ to obtain $W^\delta(t) \in G$. Here $W^\delta(0-) = X(0)$. The size of the jump is such that $W^\delta(t)$ is a strictly positive distance (depending on $\delta$) from the boundary of $G$. The jump vector is obtained as a measurable function of $W^\delta(t-)$. To ensure the measurability, each point $x$ on $\partial G$ is associated with a nearby point $\bar{x}$, chosen in a measurable way from a fixed countable set of points in $\partial G$. The jump vector for $x$ is one associated with $\bar{x}$. We now specify the mapping $x \to \bar{x}$ and the associated jump vector more precisely.

By assumption (A5)(ii), for each $x \in \partial G$, there is $c(x) \in \mathbb{R}_+^{\mathbf{I}}$ such that

$$(73) \qquad \sum_{i \in \mathcal{I}(x)} c_i(x) = 1 \quad \text{and} \quad \min_{j \in \mathcal{I}(x)} \left\langle \sum_{i \in \mathcal{I}(x)} c_i(x)\gamma^i(x), n^j(x) \right\rangle \geq a.$$

By (73), Lemma 2.1 and the fact that $n^i(\cdot)$ is continuous on $\partial G_i$ for each $i \in \mathcal{I}$, we have that for each $x \in \partial G$ there is $r_x \in (0, \delta)$ such that for each $y \in B_{r_x}(x) \cap \partial G$,

$$(74) \qquad \mathcal{I}(y) \subset \mathcal{I}(x)$$

and

$$(75) \qquad \min_{j \in \mathcal{I}(x)} \left\langle \sum_{i \in \mathcal{I}(x)} c_i(x)\gamma^i(x), n^j(y) \right\rangle \geq \frac{a}{2}.$$

It follows, using the $C^1$ nature of $\partial G_i$ and the fact that $n^i(y)$ is the inward unit normal to $\partial G_i$ at $y \in \partial G$ for each $i \in \mathcal{I}(y)$, that (by choosing $r_x$ even smaller if necessary) for each $x \in \partial G$ there is $m(x) > 0$ and $r_x \in (0, \delta)$ such that for each $y \in B_{r_x}(x) \cap \partial G$, (74)–(75) hold and

$$(76) \qquad y + \lambda \sum_{i \in \mathcal{I}(x)} c_i(x)\gamma^i(x) \in G \qquad \text{for all } \lambda \in (0, m(x)).$$

Let $B^o_{r_x}(x)$ denote the interior of the closed ball $B_{r_x}(x)$ for each $x \in \partial G$. The collection $\{B^o_{r_x}(x) : x \in \partial G\}$ is an open cover of $\partial G$ and it follows that there is a countable set $\{x_k\}$ such that $\partial G \subset \bigcup_k B_{r_{x_k}}(x_k)$ and $\{x_k\} \cap B_N(0)$ is a finite set for each integer $N \geq 1$. We can further choose the set $\{x_k\}$ to be minimal in the sense that for each strict subset $C$ of $\{x_k\}$, $\{B_{r_x}(x) : x \in C\}$ does not cover $\partial G$. Let $D_k = (B_{r_{x_k}}(x_k) \setminus (\bigcup_{i=1}^{k-1} B_{r_{x_i}}(x_i)) \cap \partial G$ for each $k$. Then $D_k \neq \varnothing$ for each $k$, $\{D_k\}$ is a partition of $\partial G$, and for each $x \in \partial G$ there is a unique index $i(x)$ such that $x \in D_{i(x)}$. For each $x \in \mathbb{R}^d$, let

$$\bar{x} = \begin{cases} x, & \text{if } x \notin \partial G, \\ x_{i(x)}, & \text{if } x \in \partial G. \end{cases}$$

Note that for all $x \in \mathbb{R}^d$,

$$(77) \qquad \|x - \bar{x}\| < \delta.$$

For each $i \in \mathcal{I}$ and $x \in \mathbb{R}^d$, let

$$(78) \qquad \gamma^{i,\delta}(x) = \gamma^i(\bar{x}).$$

The mapping $x \to \bar{x}$ is Borel measurable on $\mathbb{R}^d$ and hence $\gamma^{i,\delta}$ is a Borel measurable function from $\mathbb{R}^d$ into $\mathbb{R}^d$.

We construct $(W^\delta, Y^\delta)$ as follows. Let $X$ defined on some filtered probability space $(\Omega, \mathcal{F}, \{\mathcal{F}_t\}, P)$ be a $d$-dimensional $\{\mathcal{F}_t\}$-Brownian motion with drift $\mu$ and covariance matrix $\Gamma$ such that $X$ is continuous surely and $X(0)$ has distribution $\nu$. Let

$$\tau_1 = \inf\{t \geq 0 : X(t) \in \partial G\}$$

and

$$W^\delta(t) = X(t), \qquad Y^\delta(t) = 0 \qquad \text{for } 0 \leq t < \tau_1.$$

Note that $W^\delta(\tau_1-)$ exists on $\{\tau_1 < \infty\}$ since $X$ has continuous paths and in the case that $\tau_1 = 0$, $W^\delta(0-) \equiv X(0)$. On $\{\tau_1 < \infty\}$, define

$$Y^\delta_i(\tau_1) = \begin{cases} 0, & i \notin \mathcal{I}(\overline{W^\delta(\tau_1-)}), \\ c_i(\overline{W^\delta(\tau_1-)})\left(\dfrac{m(\overline{W^\delta(\tau_1-)})}{2} \wedge \delta\right), & i \in \mathcal{I}(\overline{W^\delta(\tau_1-)}), \end{cases}$$

and

$$W^\delta(\tau_1) = X(\tau_1)$$
$$+ \left(\frac{m(\overline{W^\delta(\tau_1-)})}{2} \wedge \delta\right)\left(\sum_{i \in \mathcal{I}(\overline{W^\delta(\tau_1-)})} c_i(\overline{W^\delta(\tau_1-)})\gamma^{i,\delta}(W^\delta(\tau_1-))\right).$$

So $W^\delta, Y^\delta$ have been defined on $[0, \tau_1)$ and at $\tau_1$ on $\{\tau_1 < \infty\}$, such that:

(i) $W^\delta(t) = X(t) + \sum_{i \in \mathcal{I}} \gamma^{i,\delta}(W^\delta(0-))Y_i^\delta(0) + \sum_{i \in \mathcal{I}} \int_{(0,t]} \gamma^{i,\delta}(W^\delta(s-))\, dY_i^\delta(s)$ for all $t \in [0, \tau_1] \cap [0, \infty)$, where $W^\delta(0-) = X(0)$,

(ii) $W^\delta(t) \in \overline{G}$ for $t \in [0, \tau_1] \cap [0, \infty)$,

(iii) for $i \in \mathcal{I}$,

  (a) $Y_i^\delta(0) \geq 0$,
  (b) $Y_i^\delta$ is nondecreasing on $[0, \tau_1] \cap [0, \infty)$,
  (c) $Y_i^\delta(t) = Y_i^\delta(0) + \int_{(0,t]} 1_{\{W^\delta(s) \in U_{2\delta}(\partial G_i \cap \partial G)\}}\, dY_i^\delta(s)$ for $t \in [0, \tau_1] \cap [0, \infty)$,

(iv) $\|\Delta Y^\delta(t)\| \equiv \|Y^\delta(t) - Y^\delta(t-)\| \leq \delta$ for $t \in [0, \tau_1] \cap [0, \infty)$, where $Y^\delta(0-) \equiv 0$.

Note that (iii)(c) above contains the expression $W^\delta(s) \in U_{2\delta}(\partial G_i \cap \partial G)$. The reader may wonder why $2\delta$ appears instead of $\delta$. The reason is that at a jump time $s$ of $Y_i^\delta$, $\overline{W^\delta(s-)} \in \partial G_i \cap \partial G$ and so

$$\operatorname{dist}(W^\delta(s), \partial G_i \cap \partial G) \leq \|W^\delta(s) - W^\delta(s-)\| + \|W^\delta(s-) - \overline{W^\delta(s-)}\|$$
$$\leq \delta + \delta = 2\delta.$$

Proceeding by induction, we assume that for some $n \geq 2$, $\tau_1 \leq \cdots \leq \tau_{n-1}$ have been defined, and $W^\delta, Y^\delta$ have been defined on $[0, \tau_{n-1})$ and at $\tau_{n-1}$ on $\{\tau_{n-1} < \infty\}$, such that (i)–(iv) above hold with $\tau_{n-1}$ in place of $\tau_1$. Then we define $\tau_n = \infty$ on $\{\tau_{n-1} = \infty\}$, and on $\{\tau_{n-1} < \infty\}$ we define

$$\tau_n = \inf\{t \geq \tau_{n-1} : W^\delta(\tau_{n-1}) + X(t) - X(\tau_{n-1}) \in \partial G\}.$$

For $\tau_{n-1} \leq t < \tau_n$, let

$$Y^\delta(t) = Y^\delta(\tau_{n-1}),$$
$$W^\delta(t) = W^\delta(\tau_{n-1}) + X(t) - X(\tau_{n-1}),$$

and on $\{\tau_n < \infty\}$, let

$$Y_i^\delta(\tau_n) = \begin{cases} Y_i^\delta(\tau_{n-1}), & i \notin \mathcal{I}(\overline{W^\delta(\tau_n-)}), \\ Y_i^\delta(\tau_{n-1}) + c_i(\overline{W^\delta(\tau_n-)})\left(\dfrac{m(\overline{W^\delta(\tau_n-)})}{2} \wedge \delta\right), & i \in \mathcal{I}(\overline{W^\delta(\tau_n-)}), \end{cases}$$

and

$$W^\delta(\tau_n) = W^\delta(\tau_n-)$$
$$+ \left(\frac{m(\overline{W^\delta(\tau_n-)})}{2} \wedge \delta\right)\left(\sum_{i \in \mathcal{I}(\overline{W^\delta(\tau_n-)})} c_i(\overline{W^\delta(\tau_n-)})\gamma^{i,\delta}(W^\delta(\tau_n-))\right).$$

In this way, $W^\delta, Y^\delta$ have been defined on $[0, \tau_n)$ and at $\tau_n$ on $\{\tau_n < \infty\}$ such that (i)–(iv) hold with $\tau_n$ in place of $\tau_1$.

By construction $\{\tau_n\}_{n=1}^\infty$ is a nondecreasing sequence of stopping times. Let $\tau = \lim_{n \to \infty} \tau_n$. On $\{\tau = \infty\}$, the construction of $(W^\delta, Y^\delta)$ is complete. We now show that $\{\tau < \infty\} = \varnothing$. In fact, if $\{\tau < \infty\} \neq \varnothing$, let $\omega \in \{\tau < \infty\}$. The above construction gives $(W^\delta(\cdot, \omega), Y^\delta(\cdot, \omega))$ on the time interval $[0, \tau(\omega))$. For each $t \in [0, \tau(\omega))$, we have

$$W^\delta(t, \omega) = X(t, \omega) + \sum_{i \in \mathcal{I}} \gamma^{i,\delta}(W^\delta(0-, \omega)) Y_i^\delta(0, \omega)$$
(79)
$$+ \sum_{i \in \mathcal{I}} \int_{(0,t]} \gamma^{i,\delta}(W^\delta(s-, \omega)) \, dY_i^\delta(s, \omega).$$

Since $X$ is continuous on $[0, \infty)$, $\|\gamma^{i,\delta}(x)\| = 1$ for each $x \in \mathbb{R}^d$ and $\sum_{i \in \mathcal{I}} Y_i^\delta(0, \omega) \leq \delta$, there are constants $\tilde{\lambda} \in (0, \tau(\omega))$ and $\widetilde{M} > 0$ (depending on $\omega$) such that

(80) $$w_{\tau(\omega)}(X(\cdot, \omega) + \gamma^{i,\delta}(W^\delta(0-, \omega)) Y_i^\delta(0, \omega), \tilde{\lambda}) < \tilde{\varepsilon}$$

and

(81) $$\sup_{0 \leq t \leq \tau(\omega)} \left\| X(\cdot, \omega) + \sum_{i \in \mathcal{I}} \gamma^{i,\delta}(W^\delta(0-, \omega)) Y_i^\delta(0, \omega) \right\| \leq \widetilde{M},$$

where $w.(\cdot, \cdot)$ is defined in (3). By the choice of $\tilde{\varepsilon}, \delta$ made at the beginning of this proof, (77)–(78) and the uniform Lipschitz property of the $\gamma^i(\cdot), i \in \mathcal{I}$, it follows that (39) and (43) hold with $\gamma^{i,\delta}(y)$ and $2\delta$ in place of $\gamma^{i,n}(y, x)$ and $\delta^n$, respectively. Then by similar pathwise analysis to that used in Case 1 and 2 of the proof of Theorem 4.2, with $\widetilde{W}^n = W^n = W^\delta$, $\alpha^n = 0$, $\gamma^{i,n}(y, x) = \gamma^{i,\delta}(y)$ for each $i \in \mathcal{I}$ and $x, y \in \mathbb{R}^d$, $X^n = X + \sum_{i \in \mathcal{I}} \gamma^{i,\delta}(W^\delta(0-)) Y_i^\delta(0)$, $Y^n = Y^\delta$, $\widetilde{Y}^n = Y^\delta - Y^\delta(0)$, $\beta^n = Y^\delta(0)$ and $\delta^n = 2\delta$, we obtain that (71) holds for any $T < \tau(\omega)$ with $\omega^n = \omega$, $N_{\eta,T} = ([\tau(\omega)/\tilde{\lambda}] + 1)\varepsilon + \widetilde{M}$. It follows that $\sup_{i \in \mathcal{I}} \sup_{s \in [0,\tau(\omega))} Y_i^\delta(s, \omega)$ is finite. By the nondecreasing property of $Y_i^\delta(\cdot, \omega)$ on $[0, \tau(\omega))$ for each $i \in \mathcal{I}$, $Y_i^\delta(\tau(\omega)-, \omega)$ exists and is finite for each $i \in \mathcal{I}$. Then by (79) and the continuity of $X$, we see that $W^\delta(\tau(\omega)-, \omega)$ exists and is finite. By the construction of $Y^\delta$ and the fact that $\sum_{i \in \mathcal{I}(x)} c_i(x) = 1$ for all $x \in \partial G$, we have that

(82) $$\sum_{i \in \mathcal{I}} Y_i^\delta(\tau(\omega)-, \omega) = \sum_{n=1}^\infty \frac{m(\overline{W^\delta(\tau_n(\omega)-, \omega)})}{2} \wedge \delta.$$

Since $\tau_n(\omega) \uparrow \tau(\omega)$ as $n \to \infty$ and $W^\delta(\tau(\omega)-, \omega)$ exists, it follows that $\{W^\delta(\tau_n(\omega)-, \omega)\}_{n=1}^\infty$ converges to $W^\delta(\tau(\omega)-, \omega) \in \partial G$ as $n \to \infty$. Consequently, $\{W^\delta(\tau_n(\omega)-, \omega)\}_{n=1}^\infty$ is a bounded sequence in $\partial G$ and so by the definition of the sets $\{D_k\}$ which form a partition of $\partial G$, there is a finite set $C$ such that

$$\{W^\delta(\tau_n(\omega)-,\omega)\}_{n=1}^\infty \subset \bigcup_{k\in C} D_k.$$

Hence,

(83) $$\inf_{n\geq 1} m(\overline{W^\delta(\tau_n(\omega)-,\omega)}) \leq \inf_{k\in C} m(x_k) > 0,$$

and so the right-hand side of (82) is infinite. On the other hand, since $\sup_{i\in\mathcal{I}}\sup_{s\in[0,\tau(\omega))} Y_i^\delta(s,\omega)$ is finite, the left-hand side of (82) is finite. This yields the desired contradiction and so $\{\tau < \infty\} = \varnothing$ and we have constructed $(W^\delta, Y^\delta)$ on $[0,\infty)$.

From the construction above, we can see that $W^\delta$ and $Y^\delta$ are well-defined stochastic processes with sample paths in $D([0,\infty),\mathbb{R}^d)$ and $D([0,\infty),\mathbb{R}^{\mathbf{I}})$. They are adapted to the filtration generated by $X$ and satisfy (i)–(iv) above with $[0,\infty)$ in place of $[0,\tau_1]$.

Consider a sequence of sufficiently small $\delta$'s, denoted by $\{\delta^n\}$, such that $\delta^n \downarrow 0$ as $n \to \infty$. For each $\delta^n$, let $(W^{\delta^n}, Y^{\delta^n})$ be the pair constructed as above for the same process $X$. By the above properties and the fact that for each $i \in \mathcal{I}$ and $x, y \in \mathbb{R}^d$,

$$\|\gamma^{i,\delta^n}(y) - \gamma^i(x)\| \leq \|\gamma^i(\bar{y}) - \gamma^i(x)\| \leq L\|\bar{y} - x\| \leq L(\delta^n + \|y - x\|),$$

we obtain that Assumption 4.1 holds with $\widetilde{W}^n = W^n = W^{\delta^n}$, $\alpha^n = 0$, $\gamma^{i,n}(y,x) = \gamma^{i,\delta^n}(y)$ for each $i \in \mathcal{I}$ and $x, y \in \mathbb{R}^d$, $X^n = X + \sum_{i\in\mathcal{I}} \gamma^{i,\delta^n}(W^{\delta^n}(0-))Y_i^{\delta^n}(0)$, $Y^n = Y^{\delta^n}$, $\widetilde{Y}^n = Y^{\delta^n} - Y^{\delta^n}(0)$, $\beta^n = Y^{\delta^n}(0)$ and $2\delta^n$ in place of $\delta^n$. By invoking the first part of Theorem 4.3, we obtain that $\{\mathcal{Z}^{\delta^n}\}_{n=1}^\infty = \{(W^{\delta^n}, X^{\delta^n}, Y^{\delta^n})\}_{n=1}^\infty$ is $C$-tight and any weak limit point $\mathcal{Z}$ of this sequence satisfies conditions (i), (ii) and (iv) of Definition 2.1 with $\mathcal{F}_t = \sigma\{\mathcal{Z}(s): 0 \leq s \leq t\}$, $t \geq 0$. Note that condition (vi)$'$ of Theorem 4.3 holds trivially. Furthermore, $M^{\delta^n} = \{X^{\delta^n}(t) - X^{\delta^n}(0) - \mu t, t \geq 0\} = \{X(t) - X(0) - \mu t, t \geq 0\}$ is a martingale with respect to the filtration generated by $X$. Since $W^{\delta^n}$, $Y^{\delta^n}$ are adapted to this filtration, it follows that $M^{\delta^n}$ is a martingale with respect to the filtration generated by $W^{\delta^n}$, $X^{\delta^n}$, $Y^{\delta^n}$ (which in fact is the same as that generated by $X$). For each $t \geq 0$, $X^{\delta^n}(t) - X^{\delta^n}(0) = X(t) - X(0)$ and so trivially this forms a uniformly integrable sequence as $n$ varies. It follows from Proposition 4.1 that condition (vii) of Theorem 4.3 holds. Hence, any weak limit point of $\{\mathcal{Z}^{\delta^n}\}_{n=1}^\infty$ is an extended SRBM with the data $(G, \mu, \Gamma, \{\gamma^i, i \in \mathcal{I}\}, \nu)$. $\square$

5.2. *SRBMs in convex polyhedrons with constant reflection fields.* Existence and uniqueness in law for SRBMs living in convex polyhedrons with a constant reflection field on each boundary face has been studied by Dai and Williams [4]. In this subsection, we state a consequence of our invariance principle using the results in [4] to establish uniqueness in law. In this case, $\overline{G}$ is defined in terms of $\mathbf{I}$ ($\mathbf{I} \geq 1$) $d$-dimensional unit vectors $\{n^i, i \in \mathcal{I}\}$ and an $\mathbf{I}$-dimensional vector $\beta = (\beta_1, \ldots, \beta_{\mathbf{I}})'$ such that

(84) $$\overline{G} \equiv \{x \in \mathbb{R}^d : \langle n^i, x \rangle \geq \beta_i \text{ for all } i \in \mathcal{I}\}.$$

It is assumed that $G$ is nonempty and that the set $\{(n^1, \beta_1), \ldots, (n^{\mathbf{I}}, \beta_{\mathbf{I}})\}$ is minimal in the sense that no proper subset defines $\overline{G}$. For each $i \in \mathcal{I}$, let $F_i$ denote the boundary face: $\{x \in \overline{G} : \langle n^i, x \rangle = \beta_i\}$. Then, $n^i$ is the inward unit normal to $F_i$. A constant vector field $\gamma^i$ of unit length specifies the direction of reflection associated with $F_i$.

DEFINITION 5.2. For each $\varnothing \neq \mathcal{K} \subset \mathcal{I}$, define $F_{\mathcal{K}} = \bigcap_{i \in \mathcal{K}} F_i$. Let $F_{\varnothing} = \overline{G}$. A set $\mathcal{K} \subset \mathcal{I}$ is maximal if $\mathcal{K} \neq \varnothing$, $F_{\mathcal{K}} \neq \varnothing$ and $F_{\mathcal{K}} \neq F_{\bar{\mathcal{K}}}$ for any $\bar{\mathcal{K}} \supset \mathcal{K}$ such that $\bar{\mathcal{K}} \neq \mathcal{K}$.

In [4], Dai and Williams introduced the following assumption.

ASSUMPTION 5.1. *For each maximal $\mathcal{K} \subset \mathcal{I}$,*

(S.a) *there is a positive linear combination $n = \sum_{i \in \mathcal{K}} b_i n^i$ ($b_i > 0\ \forall i \in \mathcal{K}$) of the $\{n^i, i \in \mathcal{K}\}$ such that $\langle n, \gamma^i \rangle > 0$ for all $i \in \mathcal{K}$,*
(S.b) *there is a positive linear combination $\gamma = \sum_{i \in \mathcal{K}} c_i \gamma^i$ ($c_i > 0\ \forall i \in \mathcal{K}$) of the $\{\gamma^i, i \in \mathcal{K}\}$ such that $\langle n^i, \gamma \rangle > 0$ for all $i \in \mathcal{K}$.*

REMARK. For the given $\overline{G}$ and constant vector fields $\{\gamma^i, i \in \mathcal{I}\}$, Assumption 5.1 is equivalent to assumption (A5).

DEFINITION 5.3. The convex polyhedron $\overline{G}$ is simple if for each $\mathcal{K} \subset \mathcal{I}$ such that $\mathcal{K} \neq \varnothing$ and $F_{\mathcal{K}} \neq \varnothing$, exactly $|\mathcal{K}|$ distinct faces contain $F_{\mathcal{K}}$.

REMARK. The polyhedron $\overline{G}$ is simple if and only if $\mathcal{K}$ is maximal for every $\mathcal{K}$ such that $\varnothing \neq \mathcal{K} \subset \mathcal{I}$ and $F_{\mathcal{K}} \neq \varnothing$. It is shown in [4] that when $\overline{G}$ is simple, (S.a) holds for all maximal $\mathcal{K} \subset \mathcal{I}$ if and only if (S.b) holds for all maximal $\mathcal{K} \subset \mathcal{I}$.

Dai and Williams [4] showed that Assumption 5.1 is sufficient for existence and uniqueness in law of SRBMs living in $\overline{G}$ with the reflection fields $\{\gamma^i, i \in \mathcal{I}\}$ and fixed starting point. [They also showed that condition (S.b) holding for all maximal $\mathcal{K} \subset \mathcal{I}$ is necessary for existence of an SRBM starting from each point in $\overline{G}$. Consequently, when $\overline{G}$ is simple, Assumption 5.1 is necessary and sufficient for existence of an SRBM starting from each point in $\overline{G}$.] This yields the following consequence of our invariance principle.

THEOREM 5.4. *Let $G$ be a nonempty domain such that $\overline{G}$ is a convex polyhedron of the form* (84) *(with minimal description), and let $\{\gamma^i, i \in \mathcal{I}\}$ be a family of constant vector fields of unit length satisfying Assumption* 5.1. *Suppose that Assumption* 4.1 *and* (vi)′, (vii) *of Theorem* 4.3 *hold. Then $W^n \Rightarrow W$ as $n \to \infty$ where $W$ is an SRBM associated with $(G, \mu, \Gamma, \{\gamma^i, i \in \mathcal{I}\}, \nu)$.*

PROOF. Clearly (A1) holds. Assumptions (A2)–(A3) hold by Lemma A.3. Since for each $i \in \mathcal{I}$, $\gamma^i(\cdot)$ is a constant vector field of unit length, assumption (A4) holds trivially. Assumption (A5) is implied by Assumption 5.1. Hence by Theorem 4.3, the only thing that we have to check is condition (viii) of Theorem 4.3, that is, uniqueness in law for SRBMs in convex polyhedrons with constant reflection fields of unit length. But this is proved in Theorem 1.3 of [4] for a fixed starting point in $\overline{G}$ and follows by a standard conditioning argument for the initial distribution $\nu$. □

5.3. *SRBMs in bounded domains with piecewise smooth boundaries.* Dupuis and Ishii [6] have established sufficient conditions for the existence and pathwise uniqueness of reflecting diffusions living in the closures of bounded domains with piecewise smooth boundaries. In this subsection, we state a consequence of our invariance principle using the results in [6] to establish uniqueness in law.

THEOREM 5.5. *Let $G$ be a bounded domain and $\{\gamma^i, i \in \mathcal{I}\}$ be a family of reflection fields that satisfy assumptions* (A1)–(A4) *and* (A5)′ *in Section* 3. *We further assume that for each $i \in \mathcal{I}$, $\gamma^i(\cdot)$ is once continuously differentiable with locally Lipschitz continuous first partial derivatives. Suppose that Assumption* 4.1 *and* (vi)′, (vii) *of Theorem* 4.3 *hold. Then $W^n \Rightarrow W$ as $n \to \infty$ where $W$ is an SRBM associated with $(G, \mu, \Gamma, \{\gamma^i, i \in \mathcal{I}\}, \nu)$.*

REMARK. We remind the reader that in view of Lemma 3.1, to verify condition (A5)′, one only needs to show that (i) or (ii) holds for all $x \in \partial G$. However, as can be seen from the proof below, both forms of the condition can be useful.

PROOF OF THEOREM 5.5. This theorem follows from Theorem 4.3 and uniqueness in law for the associated SRBMs. The latter follows by a standard argument from the pathwise uniqueness established in Corollary 5.2 of [6] for their Case 2. The conditions required for that case are satisfied in particular because (A5)′(ii) implies condition (3.8) of [6]. That condition (3.8) readily implies condition (3.6) of [6]; and, by [5], under the additional smoothness assumptions imposed on the $\gamma^i$ in the statement of our theorem, condition (3.8) also implies condition (3.7) in [6]. In addition, (A5)′(i) implies that for each $x \in \partial G$, $\langle \gamma^i(x), n^i(x) \rangle > 0$ for each $i \in \mathcal{I}(x)$, and furthermore, since (A5)′ implies (A5), we have by (A5)(i) that the origin does not belong to the convex hull of the $\{\gamma^i(x) : i \in \mathcal{I}(x)\}$. □

## APPENDIX: AUXILIARY LEMMAS

LEMMA A.1. *Suppose that $G$ is bounded. If assumption* (A1) *holds, then assumption* (A2) *holds.*

PROOF. To see this, suppose $G$ is bounded and assumption (A1) holds. Fix $\varepsilon \in (0,1)$. For each $i \in \mathcal{I}$ and $z \in \partial G_i \cap \partial G$, by the $C^1$ property of $\partial G_i$, there is a neighborhood $V_z$ of $z$ and a constant $R(\varepsilon, i, z) > 0$ such that for all $x \in V_z \cap \partial G_i \cap \partial G$ and $y \in \overline{G}_i$ such that $\|x - y\| < R(\varepsilon, i, z)$,

$$\langle n^i(x), y - x \rangle \geq -\varepsilon \|y - x\|. \tag{85}$$

Assumption (A2) then follows by a standard compactness argument. □

LEMMA A.2. *Suppose that $G$ is a nonempty bounded domain satisfying* (5), *where for each $i \in \mathcal{I}$, $G_i$ is a nonempty domain. Then assumption* (A3) *holds.*

PROOF. We prove the lemma by contradiction. Suppose that assumption (A3) does not hold. Then, since there are only finite many $\mathcal{J} \subset \mathcal{I}$, $\mathcal{J} \neq \varnothing$, there is an $\varepsilon > 0$, a nonempty set $\mathcal{J} \subset \mathcal{I}$, a sequence $\{r_n\} \subset (0, \infty)$ with $r_n \to 0$ as $n \to \infty$, a sequence $\{x_n\} \subset \mathbb{R}^d$ such that for each $n$, $x_n \in \bigcap_{j \in \mathcal{J}} U_{r_n}(\partial G_j \cap \partial G)$ and $\mathrm{dist}(x_n, \bigcap_{j \in \mathcal{J}} (\partial G_j \cap \partial G)) > \varepsilon$. But since $\overline{G}$ is bounded, $\{x_n\}$ is bounded and without loss of generality we may assume that $x_n \to x$ as $n \to \infty$ for some $x \in \mathbb{R}^d$. It follows that $x \in \bigcap_{j \in \mathcal{J}} (\partial G_j \cap \partial G)$, since for each $j \in \mathcal{J}$,

$$\mathrm{dist}(x, \partial G_j \cap \partial G) \leq \|x_n - x\| + \mathrm{dist}(x_n, \partial G_j \cap \partial G) \leq \|x_n - x\| + r_n \to 0$$

as $n \to \infty$. This is inconsistent with $x_n \to x$ and $\mathrm{dist}(x_n, \bigcap_{j \in \mathcal{J}} (\partial G_j \cap \partial G)) > \varepsilon$. □

LEMMA A.3. *Suppose* (A1) *holds where*

$$G_i = \{x \in \mathbb{R}^d : \langle n^i, x \rangle > \beta_i\} \qquad \text{for } i \in \mathcal{I}, \tag{86}$$

$\{n^i, i \in \mathcal{I}\}$ *is a finite collection of $d$-dimensional vectors of unit length, and for $\mathbf{I} = |\mathcal{I}|$, $\beta = (\beta_1, \ldots, \beta_\mathbf{I})'$ is an $\mathbf{I}$-dimensional vector. (Thus, $\overline{G}$ is a convex polyhedron.) Assume that for each $i \in \mathcal{I}$, $\partial G_i \cap \partial G \neq \varnothing$. Then assumptions* (A2) *and* (A3) *hold.*

PROOF. Assumption (A2) holds automatically since $G$ is convex. In order to show that assumption (A3) holds, we just need to show that for each $\mathcal{J} \subset \mathcal{I}$ with $\mathcal{J} \neq \varnothing$,

$$(87) \quad \sup\left\{\operatorname{dist}\left(x, \bigcap_{j \in \mathcal{J}}(\partial G_j \cap \partial G)\right) : x \in \bigcap_{j \in \mathcal{J}} U_r(\partial G_j \cap \partial G)\right\} \to 0$$

as $r \to 0$. Fix $\mathcal{J} \subset \mathcal{I}$ such that $\mathcal{J} \neq \varnothing$. Then $\bigcap_{j \in \mathcal{J}}(\partial G_j \cap \partial G)$ is the collection of all solutions $x \in \mathbb{R}^d$ to the following system of linear inequalities:

(LS)
$$\langle n^i, x \rangle \geq \beta_i \quad \text{for all } i \in \mathcal{I},$$
$$\langle -n^i, x \rangle \geq -\beta_i \quad \text{for all } i \in \mathcal{J}.$$

Suppose that $\bigcap_{j \in \mathcal{J}}(\partial G_j \cap \partial G) \neq \varnothing$, that is, (LS) has at least one solution. By a theorem of Hoffman [11], with supporting lemmas proved by Agmon [1], there is a constant $C > 0$ (depending only on $\{n^i, i \in \mathcal{I}\}$ and not on $\beta$) such that for any $x \in \mathbb{R}^d$ there exists a solution $x_0 \in \mathbb{R}^d$ of (LS) with

$$(88) \quad \|x - x_0\| \leq C\left(\sum_{i \in \mathcal{I}}(\beta_i - \langle n^i, x\rangle)^+ + \sum_{i \in \mathcal{J}}(-\beta_i - \langle -n^i, x\rangle)^+\right).$$

For $r > 0$, any $x \in \bigcap_{j \in \mathcal{J}} U_r(\partial G_j \cap \partial G)$ satisfies the following:

($r$-LS)
$$\langle n^i, x \rangle \geq \beta_i - r \quad \text{for all } i \in \mathcal{I},$$
$$\langle -n^i, x \rangle \geq -\beta_i - r \quad \text{for all } i \in \mathcal{J}.$$

Then by (88), there is $x_0 \in \bigcap_{j \in \mathcal{J}}(\partial G_j \cap \partial G)$ such that

$$\operatorname{dist}\left(x, \bigcap_{j \in \mathcal{J}}(\partial G_j \cap \partial G)\right) \leq \|x - x_0\| \leq 2C|\mathcal{I}|r.$$

It follows that (87) holds when $\bigcap_{j \in \mathcal{J}}(\partial G_j \cap \partial G) \neq \varnothing$.

Now suppose that $\bigcap_{j \in \mathcal{J}}(\partial G_j \cap \partial G) = \varnothing$, that is, (LS) has no solution. We shall use an argument by contradiction to show that $\bigcap_{j \in \mathcal{J}} U_r(\partial G_j \cap \partial G) = \varnothing$ for all $r$ sufficiently small. Suppose that this is not true. Then we have that $\bigcap_{j \in \mathcal{J}} U_r(\partial G_j \cap \partial G) \neq \varnothing$ for all $r \in (0, \infty)$. As we have seen before, any $x \in \bigcap_{j \in \mathcal{J}} U_r(\partial G_j \cap \partial G)$ is a solution to ($r$-LS). We now construct a Cauchy sequence. Let $x_1 \in \bigcap_{j \in \mathcal{J}} U_{1/2}(\partial G_j \cap \partial G)$. Then $x_1$ is a solution to ($\frac{1}{2}$-LS). Since ($\frac{1}{2^2}$-LS) has at least one solution, by the theorem of Hoffman [11] (using the fact that the constant $C$ depends only on $\{n^i, i \in \mathcal{I}\}$), we conclude that there is a solution $x_2$ to ($\frac{1}{2^2}$-LS) such that $\|x_1 - x_2\| \leq \frac{C'}{2^2}$, where $C' = 2C|\mathcal{I}|$. Continuing in this manner, we can obtain a sequence $\{x_n\}_{n=1}^{\infty}$ such that for each $n \geq 1$, $\|x_n - x_{n+1}\| \leq \frac{C'}{2^{n+1}}$ and $x_{n+1}$ is a solution of ($\frac{1}{2^{n+1}}$-LS). The sequence $\{x_n\}_{n=1}^{\infty}$ is Cauchy. Hence, there is an $x^* \in \mathbb{R}^d$ such that $x_n \to x^*$ as $n \to \infty$, and $x^*$ is a solution to (LS). This contradicts the supposition that $\bigcap_{j \in \mathcal{J}} (\partial G_j \cap \partial G) = \varnothing$. Thus we have that $\bigcap_{j \in \mathcal{J}} U_r(\partial G_j \cap \partial G) = \varnothing$ for all $r$ sufficiently small, and for such $r$,

$$\sup\left\{\operatorname{dist}\left(x, \bigcap_{j \in \mathcal{J}} (\partial G_j \cap \partial G)\right) : x \in \bigcap_{j \in \mathcal{J}} U_r(\partial G_j \cap \partial G)\right\} = 0,$$

by convention.

Combining the above we see that for each $\mathcal{J} \subset \mathcal{I}$ with $\mathcal{J} \neq \varnothing$, (87) holds and hence assumption (A3) holds. $\square$

REMARK. In fact, under the assumptions of Lemma A.3, there is a constant $C > 0$ such that $D(u) \leq Cu$ for each $u \geq 0$ and $D(\cdot)$ defined as in assumption (A3).

LEMMA A.4. *Given $T > 0$, functions $\phi, \{\phi^n\}_{n=1}^{\infty}$ in $D([0,\infty), \mathbb{R}^d)$, and $\chi, \{\chi^n\}_{n=1}^{\infty}$ in $D([0,\infty), \mathbb{R})$, suppose that $\sup_{0 \leq s \leq T} \|\phi^n(s) - \phi(s)\| \to 0$ and $\sup_{0 \leq s \leq T} |\chi^n(s) - \chi(s)| \to 0$ as $n \to \infty$. Assume that $\chi^n$ is nondecreasing for each $n$. Then for any sequence of real valued continuous functions $\{f^n\}_{n=1}^{\infty}$ defined on $\mathbb{R}^d$ such that $f^n$ converges uniformly on each compact set to a continuous function $f : \mathbb{R}^d \to \mathbb{R}$, we have*

$$(89) \qquad \int_{(0,t]} f^n(\phi^n(s)) \, d\chi^n(s) \to \int_{(0,t]} f(\phi(s)) \, d\chi(s) \qquad as \ n \to \infty,$$

*uniformly for $t \in [0, T]$.*

PROOF. By replacing $\chi^n(\cdot)$ and $\chi(\cdot)$ by $\chi^n(\cdot) - \chi^n(0)$ and $\chi(\cdot) - \chi(0)$, respectively, we may assume that $\chi^n(0) = \chi(0) = 0$. It is straightforward to see by the uniform convergence of $\{\chi^n\}$ to $\chi$ on $[0, T]$ that $\chi$ inherits the nondecreasing property of the $\{\chi^n\}$.

By the triangle inequality,

$$(90) \quad \begin{aligned} \sup_{0 \leq t \leq T} &\left|\int_{(0,t]} f^n(\phi^n(s)) \, d\chi^n(s) - \int_{(0,t]} f(\phi(s)) \, d\chi(s)\right| \\ &\leq \sup_{0 \leq t \leq T} \left|\int_{(0,t]} (f^n(\phi^n(s)) - f(\phi(s))) \, d\chi^n(s)\right| \\ &\quad + \sup_{0 \leq t \leq T} \left|\int_{(0,t]} f(\phi(s)) \, d(\chi^n(s) - \chi(s))\right|. \end{aligned}$$

For the first term on the right-hand side of the above inequality, we have

$$\sup_{0 \leq t \leq T} \left|\int_{(0,t]} (f^n(\phi^n(s)) - f(\phi(s))) \, d\chi^n(s)\right| \\ \leq \sup_{0 \leq s \leq T} |f^n(\phi^n(s)) - f(\phi(s))| \chi^n(T),$$

where the right-hand side member above tends to zero as $n \to \infty$ by the uniform convergence of $\phi^n$ to $\phi$ on $[0,T]$ (which implies uniform boundedness of $\{\phi^n\}$ on $[0,T]$), the uniform convergence of $f^n$ to $f$ on compact sets, the continuity of $f$, and the convergence of $\chi^n(T)$ to $\chi(T)$. For the second term, note that since $f(\phi(\cdot)) \in D([0,\infty), \mathbb{R})$, by Theorem 3.5.6, Proposition 3.5.3 and Remark 3.5.4 of [7], there is a sequence of step functions $\{z^k\}_{k=1}^{\infty}$ of the form

$$z^k(\cdot) = \sum_{i=1}^{l_k} z^k(t_i^k) 1_{[t_i^k, t_{i+1}^k)}(\cdot), \tag{91}$$

where $1 \leq l_k < \infty$, $0 = t_1^k < t_2^k < \cdots < t_{l_k+1}^k < \infty$ and $\sup_{0 \leq s \leq T} |f(\phi(s)) - z^k(s)| \to 0$ as $k \to \infty$. Then

$$\sup_{0 \leq t \leq T} \left| \int_{(0,t]} f(\phi(s)) \, d(\chi^n(s) - \chi(s)) \right|$$

$$\leq \sup_{0 \leq t \leq T} \left| \int_{(0,t]} (f(\phi(s)) - z^k(s)) \, d(\chi^n(s) - \chi(s)) \right|$$

$$+ \sup_{0 \leq t \leq T} \left| \int_{(0,t]} z^k(s) \, d(\chi^n(s) - \chi(s)) \right|$$

$$\leq \sup_{0 \leq s \leq T} |f(\phi(s)) - z^k(s)| (\chi^n(T) + \chi(T))$$

$$+ \sup_{0 \leq t \leq T} \sum_{i=1}^{l_k} |z^k(t_i^k)| |(\chi^n - \chi)((t_{i+1}^k \wedge t)-) - (\chi^n - \chi)((t_i^k \wedge t)-)|.$$

For fixed $k$, the last term above can be made as small as we like for all $n$ sufficiently large since $\chi^n \to \chi$ uniformly on $[0,T]$. The desired result follows. □

REMARK. The proof of Lemma A.4 is a modification of the proof of the related Lemma 2.4 in [4]. The difference in assumptions is that in [4] it is assumed that $\phi^n \to \phi$ in the $J_1$-topology rather than uniformly on $[0,T]$, $\chi^n, \chi \in C([0,\infty), \mathbb{R}_+)$ rather than $\chi^n, \chi \in D([0,\infty), \mathbb{R})$, and there is a single function $f$ rather than a sequence $\{f^n\}$.

DEPARTMENT OF MATHEMATICAL SCIENCES
CARNEGIE MELLON UNIVERSITY
PITTSBURGH, PENNSYLVANIA 15213
USA
E-MAIL: weikang@andrew.cmu.edu

DEPARTMENT OF MATHEMATICS
UNIVERSITY OF CALIFORNIA AT SAN DIEGO
9500 GILMAN DRIVE
LA JOLLA, CALIFORNIA 92093
USA
E-MAIL: williams@math.ucsd.edu